\documentclass[oneside,preprint,12pt,numbers,sort&compress]{elsarticle}
 \usepackage{hyperref} %for contents
 \usepackage{amssymb}
 \usepackage{amsthm}
 \usepackage{mathrsfs}
 \usepackage{amsmath}
  \usepackage{lineno}
\usepackage[text={170mm,240mm},centering]{geometry}

\journal{}

\begin{document}
 \newtheorem{theorem}{Theorem}[section]
\newtheorem{lemma}[theorem]{Lemma}
\newtheorem{proposition}[theorem]{Proposition}
\newtheorem{remark}[theorem]{Remark}
\newtheorem{definition}[theorem]{Definition}
\newtheorem{Hypothesis}[theorem]{Hypothesis}
\newtheorem{corollary}[theorem]{Corollary}
\numberwithin{equation}{section}
\allowdisplaybreaks[4]
\begin{frontmatter}
\title{$(H^2,H^3)$-global attractors and finite-dimensionality for the calmed 3D Navier-Stokes equations}

\author[mymainaddress]{Ruofei Liu}
\ead{lrf030311@163.com}
\address[mymainaddress]{School of Mathematical Sciences, University of Jinan, Jinan 250022, P.R.  China}

\author[mymainaddress]{Hui Liu\corref{mycorrespondingauthor}}
\cortext[mycorrespondingauthor]{Corresponding author}
\ead{ss\_liuhui@ujn.edu.cn}

\author[mymainaddress1]{Chengfeng Sun}
\ead{sch200130@163.com}
\address[mymainaddress1]{School of Applied Mathematics, Nanjing University of Finance and Economics, Nanjing 210023, P.R. China}

\author[mymainaddress2]{Jie Xin}
\ead{fdxinjie@sina.com}
\address[mymainaddress2]{School of Information Engineering, Shandong Youth University of Political Science, Jinan 250103, P.R. China}

\begin{abstract}
	
In this paper, the long-time dynamics for the three-dimensional incompressible constant-density calmed Navier-Stokes equations is investigated on periodic domain. And the key to "calmed" is to constrain the larger velocity field by introducing a calming function of appropriate conditions into the nonlinear term of the equations. First, by establishing the uniform estimates of the solutions, it is proved that the system has a $(V,H^2)$-global attractor for $f\in H$, and further shows that the system also has an $(H^2,H^3)$-global attractor for $f\in V$. Then, it focuses on the long-time dynamical properties of the system, such as backwards uniqueness, Hausdorff and fractal dimensions, and determining nodes and determining modes. Thus, the global attractor theory of the three-dimensional fluid equations is enriched.
\end{abstract}

\begin{keyword}  Calmed Navier-Stokes equations, Global attractor, Dimensions, Determining modes.
\MSC[2020] 35B40, 35Q30, 37L30.
\end{keyword}
\end{frontmatter}

\section{Introduction}

Navier-Stokes equations (NSE) are the  mathematical model for the flow of the incompressible fluids, and it has important applications in many fields, for example, fluid mechanics, aerodynamics and turbulence theory (see \cite{Rob01, Robinson2016}). Because of the strong nonlinear term of the equations, the existence, uniqueness, regularity and long-time dynamical behavior of the solutions of the NSE always are important research topics in partial differential equations and infinite-dimensional dynamical system. In particular, global regularity of the three-dimensional (3D) incompressible NSE is still an important and public question in mathematics (see \cite{Fefferman2000, Robinson2020}). Different from the two-dimensional (2D) situations, the nonlinear term for the 3D equations has a more complex energy transfer structure, so the classical energy estimates are usually not enough establish the required high-order estimates(see \cite{Hou2009, Farwig2009, O2023, Lei2024}). Thus, how to establish a global well-posedness for 3D fluid equations and further study its long-time evolutionary behavior is a significant direction for the current relevant research.

However, since the global existence and uniqueness of the strong solutions of the standard 3D NSE is not generally proven, the classical semigroup theory and the results related to global attractor in 2D situations cannot be directly extended to the 3D case. Faced with this difficulty, Sell \cite{Sell1996} studied the 3D NSE from the multivalued evolution system, established a theory of the weak global attractor, and considered the situation that attractors consist of strong solutions under the additional conditions. To overcome the difficulties of the classical 3D NSE in terms of global regularity, researchers proposed numerous modified models. Currently, the common ideas for constructing them mainly include introducing the damping term, truncating the nonlinear terms and changing the time evolution structure, etc, which can be seen in \cite{Caraballo2008, Romito2009, Marin-Rubio2011, Zhao2017, Caraballo2023}. These results show that the well-posedness of 3D fluid models can be improved to some extent by adjusting the nonlinear term appropriately, and providing an effective way to develop its long-time dynamics. It should be noted that the above methods are different from the way in this paper, but their research also gives an important reference.

This paper mainly focuses on the method called "calmed" and  first proposed by Mattew Enlow \cite{Enlow2024Nonlinearity}. The model reduces the effective growth of nonlinear term by introducing the calming function in advective term and local constrains on larger velocity fields. Unlike the ways of directly filtering or mollifying derivatives, calmed NSE mainly achieves constraints on the system by modifying the nonlinear term, there is not need to adjust the boundary conditions of the original equations accordingly (see \cite{Enlow2024Nonlinearity, Enlow2024, Duan2025arXiv, Enlow2024arXiv}). Under the suitable conditions of the calming function, Enlow, Larios and Wu \cite{Enlow2024} derived the theory of global well-posedness for the calmed 3D NSE, studied the relevant energy identity and convergence properties of the solutions, and proved that the existence of global attractors for the system of rotational calmed NSE. Additionally, Tihn (\cite{Tinh2026}) considered the well-posedness, convergence, determining wavenumber and regularity criteria of solutions for equations on the calming conditions.

From the viewpoint of the dynamical system, global attractors not only describe the asymptotic state of the system after long-time evolution, but also indicates the stability and long-time dynamical structure for the system (for example, \cite{Rob01, Hale1988, Raugel1989, Temam1997, Chepyzhov2002, Chueshov2002, Miranville2008, Trelat2015, Liu2021, Liu2023, Cui2022, Cui2024, Cui2024Forward, Cui2026, Dashkovskiy2026}). And in the study of the NSE, the effect of calming function on long-time regularity and strong topological attraction by the construction of different absorbing sets and proving the asymptotic compactness in various spaces (see \cite{Enlow2024arXiv, Song2011}). Meanwhile, it is useful to understand the orbit structure of the attractors by the proof of backwards uniqueness (see \cite{Temam1997, Robinson2013}). On this basis, establishing the existence of finite Hausdorff and fractal dimensions for global attractors (see \cite{Rob01, Robinson2013, Huang2018, Kou2026}), and further finding determining nodes and determining modes (see \cite{Robinson2013, Friz2001, Toi2020, Kalantarov2009, Cheskidov2018}), which both help to show that long-time asymptotic behavior of solutions can only be dictated by finite the instant values of solutions or Fourier modes under the proper conditions.

In the above context, we investigate the three-dimensional (3D) incompressible constant-density calmed Navier-Stokes equations on $\mathbb{T}^3=[0,L]^3 \subset \mathbb{R}^3$, $L>0$,
\begin{align}\label{E:zetadian}
\left\{\begin{array}{l}
\partial_{t}u+(\zeta^{\epsilon}(u)\cdot\nabla)u+\nabla p=\nu \Delta u+f(x),
\\
\nabla\cdot u=0,
\\
u(x,0)=u_0(x).
\end{array}\right.
\end{align}
Here, $u=u(x,t)\in \mathbb{T}^3$ is the kinematic fluid velocity; $t\ge 0$ denotes time; $\nu >0$ is the kinematic viscosity; $p$ denotes the pressure; $f=f(x)$ is the external body force; $\zeta^{\epsilon}$ is a calming function, which is introduced in detail in Definition \ref{D:1.1}. For a given $u(x,0)=u_{0}(x)$ is the initial data.

Following \cite{Enlow2024}, we consider the following calming functions, which include:
\begin{align*}
\zeta^{\epsilon}(x)=
	\left\{\begin{array}{l}
		\zeta^{\epsilon}_1(x) := \frac{x}{1+\epsilon|x|},
		\\
		\zeta^{\epsilon}_2(x) := \frac{x}{1+\epsilon^2 |x|^2},
		\\
		\zeta^{\epsilon}_3(x) := q^{\epsilon}(|x|)\frac{x}{|x|},
	\end{array}\right.
\end{align*}
let $\zeta^{\epsilon}_3(0)=0$, and
\begin{align*}
	q^{\epsilon}(r)=
	\left\{\begin{array}{l}
		r, \;\;\;\;\;\;\;\;\;\;\;\;\;\;\;\;\;\;\;\;\;\;\;\;\;0\leq  r<\frac{1}{\epsilon},
		\\
		-\frac{\epsilon}{2}(r-\frac{2}{\epsilon})^2 +\frac{3}{2\epsilon}, \;\; \frac{1}{\epsilon}\leq r<\frac{2}{\epsilon},
		\\
		\frac{3}{2\epsilon}, \;\;\;\;\;\;\;\;\;\;\;\;\;\;\;\;\;\;\;\;\;\;\;\;r\geq \frac{2}{\epsilon}.
	\end{array}\right.
\end{align*}
Observe that $\zeta^{\epsilon}(x) \rightarrow x$ as $\epsilon \rightarrow 0$ for every $x\in \mathbb{T}^3$, and $\zeta^{\epsilon}(x) \in C^2$. Additionally, we impose the conditions given in the Definition \ref{D:1.1} below on $\zeta^{\epsilon}(x)$.

Specifically, the main aims and contributions of this paper can be presented by the following aspects.

\begin{itemize}
	\item \textbf{Establishing global attractors between different Sobolev spaces.} The calmed NSE exists a $(V,H^2)$-global attractor, and also has an $(H^2,H^3)$-global attractor by overcoming the coupling relationship between the calming function and the high-order derivatives of velocity field in the nonlinear term.
	\item \textbf{Backwards uniqueness.} With the help of an auxiliary function, the backwards uniqueness of calm NSE holds. Namely, if two solutions have the same state at a positive time, then they must be the same at the initial moment, too.
	\item \textbf{Hausdorff and fractal dimensions of global attractors.} Through the contraction property in the finite-dimensional complementary subspace of the projection, global attractors have finite Hausdorff and fractal dimensions in $V$ and $H^2$ respectively.
	\item \textbf{Determining nodes and determining modes.} Calmed NSE admit the determining nodes and determining modes with the high-frequency spectral properties of the Stokes operator.
\end{itemize}

This paper introduces some needed  concepts and inequalities in Section 2. In Section 3, different uniform estimates of solutions are given explicitly. In Sections 4 and 5, it is proved that the system exists $(V,H^2)$-global attractor and $(H^2,H^3)$-global attractor. And the above long-time dynamical properties are showed in Sections 6,7 and 8 successively. Additionally, an appendix is showed by Section 9.

In the following, some main concepts and results are introduced.

\subsection{Main Results}

For \eqref{E:zetadian}, the weak formulation is as follows: Given $u_0 \in H$ and $f\in H$, find $u\in L^2 (0,T;V)$ which satisfies
\begin{align}
	\left\{\begin{array}{l}
		\langle \partial_{t}u,v \rangle = \langle \nu Au, v \rangle +\langle B(\zeta^{\epsilon}(u),u),v\rangle +\langle f, v\rangle, \; \forall v \in V,
		\\
		u(x,0)=u_{0}(x),
	\end{array}\right.
\end{align}
where $A$ and $B$ have been defined and explained in detail in \eqref{A} and \eqref{B}, respectively. Moreover, it is equivalent to the following equations
\begin{align}\label{E:choudeng}
	\left\{\begin{array}{l}
		\frac{du}{dt}+ \nu Au +B(\zeta^{\epsilon}(u),u)=f,
		\\
		u(x,0)=u_{0}(x).
	\end{array}\right.
\end{align}
For notational convenience, $f$ denotes $P_{\sigma}f$.

\begin{definition}\label{D:1.1}(Definition 1.2 in \cite{Enlow2024})
	We call $\zeta^{\epsilon}:\mathbb{R}^3 \rightarrow \mathbb{R}^3$ is a calming function if the following properties hold:
	\begin{enumerate}[(i)]
		\item $\zeta^{\epsilon}$ is Lipschitz continuous with Lipschitz constant 1.
		\item For every $\epsilon>0$ fixed, $\zeta^{\epsilon}$, $D\zeta^{\epsilon}$ and $D^2\zeta^{\epsilon}$ are bounded.
		\item For every $\epsilon>0$,and for any $x\in \mathbb{R}^3$, there exists $\lambda^{\epsilon}(x)\in \mathbb{R}$ such that $\zeta^{\epsilon}(x)=\lambda^{\epsilon}(x)x$. In other words, $\zeta^{\epsilon}(x)$ is assumed to be pointwise parallel to $x$.
	\end{enumerate}
\end{definition}

\begin{lemma}\label{L:1.2}(Proposition 1.4 in \cite{Enlow2024})
	For $\epsilon >0$, the calming functions $\zeta^{\epsilon}_{i}$, $i=1,2,3$, given in Definition \ref{D:1.1} admit the following concrete bounds.
	\begin{enumerate}[(1)]
		\item For $\zeta^{\epsilon}_1$,
		\begin{align*}
			\|\zeta^{\epsilon}_1\|_{L^{\infty}}=\frac{1}{\epsilon},\;  \|D\zeta^{\epsilon}_1\|_{L^{\infty}}=1\; \textup{and} \; \|D^2\zeta^{\epsilon}_1\|_{L^{\infty}}=2\epsilon.
		\end{align*}
		\item For $\zeta^{\epsilon}_2$,
		\begin{align*}
			\|\zeta^{\epsilon}_2\|_{L^{\infty}}=\frac{1}{2\epsilon},\;  \|D\zeta^{\epsilon}_2\|_{L^{\infty}}=1 \; \textup{and} \;\|D^2\zeta^{\epsilon}_2\|_{L^{\infty}}=\frac{\epsilon}{2(2-\sqrt{2})^2}.
		\end{align*}
		\item For $\zeta^{\epsilon}_3$,
		\begin{align*}
			\|\zeta^{\epsilon}_3\|_{L^{\infty}}=\frac{3}{2\epsilon},\;  \|D\zeta^{\epsilon}_3\|_{L^{\infty}}=1\; \textup{and} \;\|D^2\zeta^{\epsilon}_3\|_{L^{\infty}}=\epsilon.
		\end{align*}
	\end{enumerate}
\end{lemma}

\begin{proof}
	All the above bounds for $\zeta^{\epsilon}_{i} \; (i=1,2,3)$ can be obtained by direct calculation (see Appendix 10.1 in \cite{Enlow2024}).
	
	Next, we denote by $D\zeta^{\epsilon}_{i}$ the first derivative of $\zeta^{\epsilon}_{i}$ and $D^2\zeta^{\epsilon}_{i}$ the second derivative of $\zeta^{\epsilon}_{i}$, then their bounds are computed as follows.
	
	For $\zeta^{\epsilon}_1$,
	\begin{align*}
		D\zeta^{\epsilon}_1 (x)=\frac{1}{(1+\epsilon |x|)^2} \; \textup{and}\;D^2\zeta^{\epsilon}_1 (x)=-\frac{2(\epsilon x +\epsilon^2 x|x|)}{|x|(1+\epsilon |x|)^4},
	\end{align*}
	then,
	\begin{align*}
		\|D\zeta^{\epsilon}_1\|_{L^{\infty}}=1 \; \textup{and} \; \|D^2\zeta^{\epsilon}_1\|_{L^{\infty}}=2\epsilon.
	\end{align*}
	
	For $\zeta^{\epsilon}_2$,
		\begin{align*}
		D\zeta^{\epsilon}_2 (x)=\frac{1-\epsilon^2|x|^2}{(1+\epsilon^2 |x|^2)^2} \; \textup{and}\; D^2\zeta^{\epsilon}_2 (x)=\frac{2\epsilon^2 x(\epsilon^2 |x|^2 -3)}{(1+\epsilon^2 |x|^2)^3},
	\end{align*}
	then,
	\begin{align*}
		\|D\zeta^{\epsilon}_2\|_{L^{\infty}}=1 \; \textup{and}\; \|D^2\zeta^{\epsilon}_2\|_{L^{\infty}}=\frac{\epsilon}{2(2-\sqrt{2})^2}.
	\end{align*}
	
		For $\zeta^{\epsilon}_3$,
		\begin{align*}
				D\zeta^{\epsilon}_{3}(x)=
			\left\{\begin{array}{l}
				1, \;\;\;\;\;\;\;\;\;\;\;\;\;\;\;\;\;\;\;\;\;\;\;\;0\leq  |x|<\frac{1}{\epsilon},
				\\
				2-\epsilon|x|, \;\;\;\;\;\;\;\;\;\;\;\;\;\; \frac{1}{\epsilon}\leq |x|<\frac{2}{\epsilon},
				\\
				0, \;\;\;\;\;\;\;\;\;\;\;\;\;\;\;\;\;\;\;\;\;\;\;\;|x|\geq \frac{2}{\epsilon},
			\end{array}\right.
		\end{align*}
	and
	\begin{align*}
		D^2\zeta^{\epsilon}_{3}(x)=
		\left\{\begin{array}{l}
			0, \;\;\;\;\;\;\;\;\;\;\;\;\;\;\;\;\;\;\;\;\;\;\;\;0\leq  |x|<\frac{1}{\epsilon},
			\\
			-\frac{\epsilon x}{|x|}, \;\;\;\;\;\;\;\;\;\;\;\;\;\;\;\;\;\;\; \frac{1}{\epsilon}\leq |x|<\frac{2}{\epsilon},
			\\
			0, \;\;\;\;\;\;\;\;\;\;\;\;\;\;\;\;\;\;\;\;\;\;\;\;|x|\geq \frac{2}{\epsilon},
		\end{array}\right.
	\end{align*}
	then,
		\begin{align*}
		\|D\zeta^{\epsilon}_3\|_{L^{\infty}}=1 \; \textup{and}\; \|D^2\zeta^{\epsilon}_3\|_{L^{\infty}}=\epsilon.
	\end{align*}
\end{proof}

Next, two main results for $(V,H^2)$-global attractor and $(H^2,H^3)$-global attractor are given by the following theorems.

\begin{theorem}\label{M:1.3}
Assume that $u_0 \in V$, $f \in H$ and $\zeta^{\epsilon}$ satisfies the conditions in Definition \ref{D:1.1}. Thus, \eqref{E:choudeng} has a $(V,H^2)$-global attractor $\mathcal{A}$, which is compact and invariant in $H^2$. Additionally, it attracts all bounded subsets of $V$ in the $H^2$-norm topology.
\end{theorem}

\begin{theorem}\label{M:1.4}
Assume that $u_0 \in H^2$, $f \in V$ and $\zeta^{\epsilon}$ satisfies the conditions in Definition \ref{D:1.1}. Hence, \eqref{E:choudeng} also exists an $(H^2,H^3)$-global attractor $\mathcal{A}$, which is compact and invariant in $H^2$. Similarly, it attracts every bounded subset of $H^2$ in the $H^3$-norm topology.
\end{theorem}

Then, the main results of the long-time dynamical properties are showed as follow.

\begin{theorem}\label{M:1.5}
Assume that $u_0 , v_0$ $\in V$, if there has $T > 0$, such that $S(T)u_0 =S(T)v_0 $, then $u_0 =v_0 $.
\end{theorem}

\begin{theorem}\label{M:1.6}
On the basis of  Theorems \ref{M:1.3} and \ref{M:1.4}, the global attractor $\mathcal{A}$ exists Hausdorff and fractal dimensions in the $V$ and $H^2$ space, respectively.
\end{theorem}

\begin{theorem}\label{M:1.7}
If there is a positive constant $\delta$, such that $d(\mathcal{N})<\delta$ holds, then $\mathcal{N}$ is a set of determining nodes for \eqref{E:choudeng} on global attractor $\mathcal{A}$.
\end{theorem}

\begin{theorem}\label{M:1.8}
There has a constant $n_0$ such that the first $n$ Fourier modes are determining for all $n>n_0$.
\end{theorem}

\section{Preliminaries}

In this section, we will introduce some symbols and preliminary concepts for later use. At the same time, the standard notations and conclusions we adopt can be found in \cite{Robinson2016, Enlow2024, Duan2025arXiv, Enlow2024Nonlinearity}.

We denote the periodic domain $\mathbb{T}^3 :=[0,L]^3$. There are constants $c_1>0$ and $c_2>0$ such that
\begin{align*}
	c_1 \|Au\|_{L^2}\leq \|u\|_{H^2}\leq c_2 \|Au\|_{L^2},
\end{align*}
where $A$ is defined by \eqref{A} (see  \cite{Dauge1989, Guermond2011}).
We also take the real $H$ inner product to be given by $\left ( \cdot,\cdot \right )$, and the $H^s$ Sobolev norm to be given by $\| \cdot \|_{H^s}$, $s\ge 1$, with the explicit forms given as follows,
\begin{align*}
	(u,v):= \sum\limits_{i=1}^3 \int_{\mathbb{T}^3}u_{i}(x)v_{i}(x)dx,
	\;\;\;
	\|u\|_{H^s}^2:= \sum\limits_{0\leq |\alpha|\leq s}\|D^{\alpha}u\|^2_{L^2}.
\end{align*}
To simplify notation, we denote $L^2(\mathbb{T}^3)^3\equiv L^2$ and $H^s(\mathbb{T}^3)^3\equiv H^s$, $s\ge 1$. And we set
\begin{align*}
	\mathcal{V}=\left\{\varphi\in [C^{\infty}(\mathbb{T}^3)]^3 : \nabla\cdot\varphi = 0 \text{ and } \int_{\mathbb{T}^3} \varphi dx = 0 \right\},
\end{align*}
then let $H$ and $V$ be the closure of $\mathcal{V}$ in $L^2$ and $H^1$.
Moreover, we denote the norm in $L^{p}(\mathbb{T}^3)$ by $\|\cdot \|_{L^p}$, and it is given by
\begin{align*}
	\|u\|_{L^p} :=\big( \int_{\mathbb{T}^3} |u(x)|^p dx \big)^{\frac{1}{p}}.
\end{align*}

Let $P_{\sigma} :L^2\rightarrow H$ be the Leray-Helmholtz orthogonal projection of $L^2$ on $H$. And let $A:\mathcal{D}(A) \subset H \rightarrow H$ be the Stokes operator defined on $\mathcal{D}(A):=H^2\cap V$, where
\begin{align}\label{A}
	A:=-P_{\sigma}\Delta,
\end{align}
where $P_{\sigma}$ and $\Delta$ commute on the periodic domain, and the operator $A$ is positive-definite, self-adjoint and with the compact inverse operator $A^{-1}$ in $H$. Applying the Hilbert-Schmidt theorem yields a complete orthonormal basis $\{\omega_{j}\}_{j=1}^{\infty}$ of $H$ consisting of eigenfunctions of the compact operator $A^{-1}$, which are also eigenfunctions of $A$. The corresponding eigenvalues $\{\lambda_{j}\}_{j=1}^{\infty}$ satisfy $0<\lambda_{1}\le \lambda_{2}\le \cdots \rightarrow \infty $ (as $j\rightarrow \infty$), and $A\omega_{j}=\lambda_{j}\omega_{j}$ for each $j\in \mathbb{N}^{+}$.
Moreover, we set $\mathcal{D}(A^{\frac{3}{2}}):= H^3\cap V$, and $\|A^{\frac{3}{2}}\cdot\|_{L^2}$ is equivalent to the norm of $H^3$.

As is well known, the Poincare inequalities are given as follows:
\begin{align*}
	&\|u\|^2_{L^2} \leq \frac{1}{\lambda_1}\|\nabla u\|^2_{L^2}, \;\; \textup{for any}\; u\in V,
	\\
	&\|\nabla u\|^2_{L^2} \leq \frac{1}{\lambda_1}\|A u\|^2_{L^2}, \;\; \textup{for any}\; u\in \mathcal{D}(A).
\end{align*}
From \cite{Enlow2024arXiv}, let us recall the Gagliardo-Nirenberg inequality as follows, assume $1\le q, r\le \infty$, and $0< \gamma<1$, for $u\in L^{q}(\mathbb{T}^3)$, such that $\partial^m u \in L^{r}(\mathbb{T}^3)$, $m>s$, then
\begin{align*}
	\|\partial^s u\|_{L^p} \leq C\|\partial^m u\|^{\gamma}_{L^r}\|u\|^{1-\gamma}_{L^q},\;\; \frac{1}{p}-\frac{s}{3}=(\frac{1}{r}-\frac{m}{3})\gamma+\frac{1}{q}(1-\gamma).
\end{align*}

The nonlinear term $B: V\times V \rightarrow V'$ is defined by
\begin{align}\label{B}
	B(u,v):=P_{\sigma}((u\cdot \nabla) v).
\end{align}
 Correspondingly, we define the $b:V\times V \times V\rightarrow \mathbb{R}$ by
\begin{align}
	b(u,v,w):=\left < B(u,v),w\right >,
\end{align}
where
\begin{align*}
	b(u,v,w)=\sum\limits_{i,j=1}^3 \int_{\mathbb{T}^3}u_i \frac{\partial v_j}{\partial x_i}w_j dx.
\end{align*}

Next, we present a lemma that is used prove the theorems of determining nodes and determining modes.
\begin{lemma}\label{L:2.1}(Lemma 36 in \cite{Robinson2013})
	Assume that $X\in L^1_{loc}(0,\infty)$ is nonnegative and if it is
	\begin{align}
		\frac{dX}{dt} +aX \leq b(t),
	\end{align}
	where $b(t)\rightarrow 0$, and for all $t\geq t_0$ such that $X(t)\leq \kappa $. Then $X(t)\rightarrow 0$ as $t\rightarrow \infty$.
\end{lemma}

Additionally, throughout this paper, $C_{\nu}$ represents a positive constant that depends on $\nu$, and its value may vary in different instances.
\section{Uniform estimates of solutions}

In this section, we derive a series of necessary uniform estimates of solutions of \eqref{E:choudeng}, which are used in the proofs below.
\begin{lemma}\label{L:3.1}
	Assume that $u_{0}\in H$, $f\in H$ and $\zeta^{\epsilon}$ satisfies the conditions in Definition \ref{D:1.1}. Thus, there are constants $\kappa_1>0$, $I_1>0$, and a time
	\begin{align*}
		t_0=\max\{0,\frac{1}{\nu\lambda_{1}}\ln{\frac{(\nu\lambda_{1})^2 \|u_0\|^{2}_{L^2}}{\|f\|^2_{L^2}}}\},
	\end{align*}
such that for all $t\geq t_0$,
	\begin{align}
		\|u(t)\|_{L^2} \le \kappa_1,
	\end{align}
	\begin{align}
		\int^{t+1}_{t}\|\nabla u(s)\|_{L^2}^{2}ds \le I_1.
	\end{align}
\end{lemma}

\begin{proof}
	Taking the inner product of \eqref{E:zetadian} with $u$ in $\mathbb{T}^3$, we have
	\begin{align*}
		\frac{1}{2}\frac{d}{dt}\|u\|^{2}_{L^2}+\nu\|\nabla u\|_{L^2}^{2}=&-\int_{\mathbb{T}^3}B(\zeta^{\epsilon}(u),u)\cdot u dx + \int_{\mathbb{T}^3}f\cdot u dx
		\nonumber\\
		\leq&\
		\|f\|_{L^2} \|u\|_{L^2}
		\nonumber\\
		\leq&\ \frac{1}{2\nu\lambda_{1}}\|f\|_{L^2}^{2}+\frac{\nu\lambda_{1}}{2}\|u\|_{L^2}^{2}
		\nonumber\\
		\leq&\ \frac{1}{2\nu\lambda_{1}}\|f\|_{L^2}^{2}+\frac{\nu}{2}\|\nabla u\|_{L^2}^{2},
	\end{align*}
	where $\int_{\mathbb{T}^3}B(\xi^{\epsilon}(u),u)\cdot u dx=0$, by the condition in Definition \ref{D:1.1}.
	So,
	\begin{align}\label{E:3.3}
		\frac{d}{dt}\|u\|^{2}_{L^2}+\nu\|\nabla u\|_{L^2}^{2}\leq\frac{1}{\nu\lambda_{1}}\|f\|_{L^2}^{2},
	\end{align}
	using Poincare inequality,
	\begin{align*}
		\frac{d}{dt}\|u\|^{2}_{L^2}+\nu\lambda_{1}\|u\|_{L^2}^{2}\leq\frac{1}{\nu\lambda_{1}}\|f\|_{L^2}^{2}.
	\end{align*}
	Then, we use Gronwall inequality to get
	\begin{align*}
		\|u(t)\|^{2}_{L^2}\leq&\ e^{-\nu\lambda_{1}t}\|u_0\|_{L^2}^{2}+\frac{1}{(\nu\lambda_{1})^2}\|f\|_{L^2}^{2}(1-e^{-\nu\lambda_{1}t})
		\nonumber\\
		\leq&\ e^{-\nu\lambda_{1}t}\|u_0\|_{L^2}^{2}+\frac{1}{(\nu\lambda_{1})^2}\|f\|_{L^2}^{2},
	\end{align*}
	and there is a time $t_{0}$, which is recommended as
	\begin{align}
		t_{0}=\max\{0,\frac{1}{\nu\lambda_{1}}\ln{\frac{(\nu\lambda_{1})^2 \|u_0\|^{2}_{L^2}}{\|f\|^2_{L^2}}}\},
	\end{align}
	such that for all $t\ge t_0$
	\begin{align*}
		\|u(t)\|^2_{L^2} \le \frac{2}{(\nu\lambda_1)^2}\|f\|^2_{L^2} :=\kappa^2_1.
	\end{align*}
	
	At the same time, we integrate \eqref{E:3.3} in $[t,t+1]$ to get
	\begin{align*}
		\nu\int^{t+1}_{t}\|\nabla u(s)\|_{L^2}^{2}ds
		\le&\
		\frac{1}{\nu\lambda_1}\|f\|_{L^2}^2 + \|u(t)\|^2_{L^2}
		\nonumber\\
		\le&\
		(\frac{1}{\nu\lambda_1}+\frac{2}{(\nu\lambda_1)^2})\|f\|^{2}_{L^2}:=\nu I_1,\; \forall t\ge t_0.
	\end{align*}
	
This completes the proof.
\end{proof}

\begin{lemma}\label{L:3.2}
Assume that $u_{0}\in H$, $f\in H$ and $\zeta^{\epsilon}$ satisfies the conditions in Definition \ref{D:1.1}. Thus, there are constants $\kappa_2>0$, $I_2>0$, and a time $t_1=t_0 +1$, such that for all $t\geq t_1$,
\begin{align}
\|\nabla u(t)\|_{L^2} \leq \kappa_2 ,
\end{align}
\begin{align}
\int^{t+1}_{t}\|Au(s)\|^2_{L^2}ds\leq I_2.
\end{align}
\end{lemma}

\begin{proof}
We take \eqref{E:choudeng} and $Au$ as the inner product and can get
\begin{align*}
	\frac{1}{2}\frac{d}{dt}\|\nabla u\|^{2}_{L^2}+\nu\|Au\|_{L^2}^{2}=&-\int_{\mathbb{T}^3}B(\zeta^{\epsilon}(u),u)\cdot Au dx + \int_{\mathbb{T}^3}f\cdot Au dx
	\nonumber\\
	\leq&\
	\int_{\mathbb{T}^3}|\zeta^{\epsilon}(u)||\nabla u||Au|dx + \int_{\mathbb{T}^3}|f||Au|dx
	\nonumber\\
	\leq&\ \frac{\nu}{2}\|Au\|_{L^2}^{2}+C_{\nu}\|\zeta^{\epsilon}\|^{2}_{L^{\infty}}\|\nabla u\|^2_{L^2} +C_\nu \|f\|^2_{L^2}.
\end{align*}
After rearranging all terms, which yields
\begin{align}\label{E:3.7}
	\frac{d}{dt}\|\nabla u\|^{2}_{L^2}+\nu\|Au\|_{L^2}^{2} \leq C_{\nu}\|\zeta^{\epsilon}\|^{2}_{L^{\infty}}\|\nabla u\|^2_{L^2} +C_\nu \|f\|^2_{L^2} .
\end{align}
Integrating the inequality over $[s,t+1]$, where $s\in (t,t+1)$ and $t\ge t_0$,
\begin{align*}
	\|\nabla u(t+1)\|^{2}_{L^2} \leq \|\nabla u(s)\|^{2}_{L^2}+ C_{\nu}\|\zeta^{\epsilon}\|^{2}_{L^{\infty}}I_1 +C_\nu \|f\|^2_{L^2}.
\end{align*}
Simultaneously, we integrate with respect to $s$ over $[t,t+1]$, where $t \ge t_0$,
\begin{align*}
	\|\nabla u(t+1)\|^{2}_{L^2} \leq& \int^{t+1}_{t}\|\nabla u(s)\|^{2}_{L^2}ds + \int^{t+1}_{t}C_{\nu}\|\zeta^{\epsilon}\|^{2}_{L^{\infty}}I_1 ds +\int^{t+1}_{t}C_\nu \|f\|^2_{L^2}ds
	\\
	\leq&\
	(1+C_{\nu}\|\zeta^{\epsilon}\|^{2}_{L^{\infty}})I_1 +C_\nu \|f\|^2_{L^2} .
\end{align*}
Then for any $t\geq t_0 +1 :=t_1$, we can get
\begin{align*}
	\|\nabla u(t)\|^{2}_{L^2} \leq (1+C_{\nu}\|\zeta^{\epsilon}\|^{2}_{L^{\infty}})I_1 +C_\nu \|f\|^2_{L^2}:= \kappa_2^2.
\end{align*}

Additionally, we integrate inequality \eqref{E:3.7} between $t$ to $t+1$, we have
\begin{align*}
	\nu\int^{t+1}_{t}\|Au(s)\|_{L^2}^{2}ds \leq& \|\nabla u(t)\|^{2}_{L^2}+ C_{\nu}\int^{t+1}_{t}\|\zeta^{\epsilon}\|^{2}_{L^{\infty}}\|\nabla u(s)\|^2_{L^2}ds +C_\nu \int^{t+1}_{t}\|f\|^2_{L^2}ds
	\\
	\leq& \kappa^2_2 +C_{\nu}\|\zeta^{\epsilon}\|^{2}_{L^{\infty}}I_1 +C_\nu \|f\|^2_{L^2} :=\nu I_2.
\end{align*}

This completes the proof.
\end{proof}

\begin{lemma}\label{L:3.3}
	Assume that $u_{0}\in H$, $f\in H$ and $\zeta^{\epsilon}$ satisfies the conditions in Definition \ref{D:1.1}. Thus, there is a constant $I_3>0$ and a time $t_1=t_0 +1$ such that for all $t\geq t_1$,
	\begin{align}
		\int^{t+1}_{t}\|u_{t}(s)\|^2_{L^2}ds\leq I_3.
	\end{align}
\end{lemma}

\begin{proof}
	We take \eqref{E:choudeng} and $u_{t}$ as the inner product on $\mathbb{T}^3$
	\begin{align*}
		\frac{\nu}{2}\frac{d}{dt}\|\nabla u\|^{2}_{L^2}+\|u_{t}\|_{L^2}^{2}=&-\int_{\mathbb{T}^3}B(\zeta^{\epsilon}(u),u)\cdot u_{t} dx + \int_{\mathbb{T}^3}f\cdot u_{t} dx
		\nonumber\\
		\leq&\
		\|\zeta^{\epsilon}(u)\|_{L^{\infty}}\|\nabla u\|_{L^2}\|u_{t}\|_{L^2} + \|f\|_{L^2}\|u_{t}\|_{L^2}
		\nonumber\\
		\leq&\ \frac{1}{2}\|u_{t}\|_{L^2}^{2}+C\|\zeta^{\epsilon}\|^{2}_{L^{\infty}}\|\nabla u\|^2_{L^2} +C \|f\|^2_{L^2}.
	\end{align*}
	Rearranging all terms and adjusting the coefficients, which yields
	\begin{align}\label{E:3.9}
		\nu\frac{d}{dt}\|\nabla u\|^{2}_{L^2}+\|u_{t}\|_{L^2}^{2} \leq C\|\zeta^{\epsilon}\|^{2}_{L^{\infty}}\|\nabla u\|^2_{L^2} +C \|f\|^2_{L^2} .
	\end{align}
	Integrating \eqref{E:3.9} over $[t,t+1]$, where $t\ge t_1$, then we can get
	\begin{align*}
		\int^{t+1}_{t}\|u_{t}(s)\|_{L^2}^{2}ds \leq& \nu\|\nabla u(t)\|^{2}_{L^2}+ C\|\zeta^{\epsilon}\|^{2}_{L^{\infty}}\int^{t+1}_{t}\|\nabla u(s)\|^2_{L^2}ds +C\|f\|^2_{L^2}
		\\
		\leq& \nu\kappa^2_2 +C\|\zeta^{\epsilon}\|^{2}_{L^{\infty}}I_1 +C \|f\|^2_{L^2} :=I_3.
	\end{align*}
	
	This completes the proof.
\end{proof}

\begin{lemma}\label{L:3.4}
	Assume that $u_{0}\in H$, $f\in H$ and $\zeta^{\epsilon}$ satisfies the conditions in Definition \ref{D:1.1}. Thus, there is a constant $\kappa_3>0$ and a time $t_2=t_1 +1$ such that for all $t\geq t_2$,
	\begin{align}
		\|u_{t}(t)\|_{L^2}\leq \kappa_3.
	\end{align}
\end{lemma}

\begin{proof}
	Differentiating \eqref{E:choudeng} with respect to $t$ yields
	\begin{align}\label{E:3.11}
		\frac{d}{dt}u_{t}+ \nu Au_t +B((D\zeta^{\epsilon}(u))u_t,u)+B(\zeta^{\epsilon}(u),u_t)=0,
	\end{align}
	then taking the inner product of \eqref{E:3.11} with $u_t$, so we have
	\begin{align*}
		\frac{1}{2}\frac{d}{dt}\|u_t\|^{2}_{L^2}+\nu\|\nabla u_{t}\|_{L^2}^{2}=&-\int_{\mathbb{T}^3}B((D\zeta^{\epsilon}(u))u_t,u)\cdot u_{t} dx
		\nonumber\\
		\leq&\
		\|D\zeta^{\epsilon}(u)\|_{L^{\infty}}\|u_{t}\|^{2}_{L^4}\|\nabla u\|_{L^2}
		\nonumber\\
		\leq&\
		C\|D\zeta^{\epsilon}\|_{L^{\infty}}\|\nabla u_t\|^{\frac{3}{2}}_{L^2}\|u_{t}\|^{\frac{1}{2}}_{L^2}\|\nabla u\|_{L^2}
		\nonumber\\
		\leq&\ \frac{\nu}{2}\|\nabla u_{t}\|_{L^2}^{2}+C_{\nu}\|D\zeta^{\epsilon}\|^{4}_{L^{\infty}}\| u_t\|^2_{L^2}\|\nabla u\|^{4}_{L^2}.
	\end{align*}
	So,
	\begin{align}\label{E:3.12}
		\frac{d}{dt}\|u_t\|^{2}_{L^2}+\nu\|\nabla u_{t}\|_{L^2}^{2} \leq C_{\nu}\|D\zeta^{\epsilon}\|^{4}_{L^{\infty}}\|\nabla u\|^{4}_{L^2}\|u_t\|^2_{L^2}.
	\end{align}
	
	Integrating \eqref{E:3.12} over $[s,t+1]$, where $s\in(t,t+1)$, then we can get
	\begin{align*}
		\|u_t(t+1)\|^2_{L^2} \leq \|u_t (s)\|^2_{L^2}+C_{\nu}\int^{t+1}_{s}\|D\zeta^{\epsilon}\|^4_{L^{\infty}}\|\nabla u(\tau)\|^4_{L^2}\|u_t(\tau)\|^2_{L^2}d\tau.
	\end{align*}
	Integrating the above inequality with respect to $s$ over $[t,t+1]$ again, where $t\ge t_1$,
	\begin{align*}
		\|u_t(t+1)\|^2_{L^2}\leq&\ \int^{t+1}_{t}(1+C_{\nu}\|D\zeta^{\epsilon}\|^4_{L^{\infty}}\|\nabla u(s)\|^4_{L^2})\|u_t(s)\|^2_{L^2}ds
		\nonumber\\
		\leq&\
		(1+C_{\nu}\|D\zeta^{\epsilon}\|^4_{L^{\infty}}\kappa^4_2)I_3.
	\end{align*}
	Then for $\forall t\ge t_1 +1 :=t_2$,
	\begin{align*}
		\|u_t(t)\|^2_{L^2} \leq (1+C_\nu\|D\zeta^{\epsilon}\|^4_{L^{\infty}}\kappa^4_2)I_3:=\kappa^2_3.
	\end{align*}
	
	This completes the proof.
\end{proof}

\begin{lemma}\label{L:3.5}
	Assume that $u_{0}\in H$, $f\in H$ and $\zeta^{\epsilon}$ satisfies the conditions in Definition \ref{D:1.1}. Thus, there is a constant $\kappa_4>0$ and a time $t_2=t_1 +1$ such that for all $t\geq t_2$,
	\begin{align}
		\|Au(t)\|_{L^2}\leq \kappa_4.
	\end{align}
\end{lemma}

\begin{proof}
	Applying Minkowski's inequality, Definition \ref{D:1.1}, and Lemmas \ref{L:3.2} and \ref{L:3.4} to \eqref{E:choudeng}, where $t\ge t_2$,
	\begin{align*}
		\nu\|Au(t)\|_{L^2}\leq&\|f\|_{L^2}+\|u_t\|_{L^2}+\|B(\zeta^{\epsilon}(u),u)\|_{L^2}
		\nonumber\\
		\leq&\
		\|f\|_{L^2}+\|u_t\|_{L^2}+\|\zeta^{\epsilon}(u)\|_{L^{\infty}}\|\nabla u\|_{L^2}
		\nonumber\\
		\leq&\
		\|f\|_{L^2}+\kappa_3 +\|\zeta^{\epsilon}\|_{L^{\infty}}\kappa_2 :=\nu \kappa_4 .
	\end{align*}

	This completes the proof.
\end{proof}

\begin{lemma}\label{L:3.6}
	Assume that $u_{0}\in H$, $f\in H$ and $\zeta^{\epsilon}$ satisfies the conditions in Definition \ref{D:1.1}. Thus, there are constants $\kappa_5>0$, $I_4>0$, and a time $t_3=t_2 +1$ such that for all $t\geq t_3$,
	\begin{align}
		\|\nabla u_{t}(t)\|_{L^2}\leq \kappa_5,
	\end{align}
	\begin{align}
		\int^{t+1}_{t}\|\nabla u_{t}(s)\|_{L^2}^{2}ds\leq I_4 .
	\end{align}
\end{lemma}

\begin{proof}
	Integrating  \eqref{E:3.12} from $t$ to $t+1$, then using Definition \ref{D:1.1} and Lemmas \ref{L:3.2} and \ref{L:3.4}, which yields,
	\begin{align*}
		\nu \int^{t+1}_{t}\|\nabla u_{t}(s)\|_{L^2}^{2}ds\leq&\|u_{t}(t)\|_{L^2}^{2}+\int^{t+1}_{t}C_{\nu}\|D\zeta^{\epsilon}\|^{4}_{L^{\infty}}\|\nabla u(s)\|_{L^2}^{4}\|u_{t}(s)\|_{L^2}^{2}ds
		\nonumber\\
		\leq&\
		\kappa_3^2 +C_{\nu}\|D\zeta^{\epsilon}\|^{4}_{L^{\infty}}\kappa^2_3 \kappa^4_2 ,
	\end{align*}
	so,
	\begin{align*}
		\int^{t+1}_{t}\|\nabla u_{t}(s)\|_{L^2}^{2}ds\leq \frac{1}{\nu}(
		\kappa_3^2 +C_{\nu}\|D\zeta^{\epsilon}\|^{4}_{L^{\infty}}\kappa^2_3 \kappa^4_2):=I_4, \; \textup{for} \; t\ge t_2.
	\end{align*}
	
	Taking the inner product of \eqref{E:3.11} with $Au_t$, where $t\ge t_2$,
	\begin{align*}
		\frac{1}{2}\frac{d}{dt}\|\nabla u_t\|^{2}_{L^2}+\nu\|Au_{t}\|_{L^2}^{2}=&-b((D\zeta^{\epsilon}(u))u_t ,u,Au_t )-b(\zeta^{\epsilon}(u) ,u_t ,Au_t )
		\nonumber\\
		\leq&\
		\|D\zeta^{\epsilon}(u)\|_{L^{\infty}}\|u_{t}\|_{L^6}\|\nabla u\|_{L^3}\|Au_t \|_{L^2}+\|\zeta^{\epsilon}(u)\|_{L^{\infty}}\|\nabla u_t \|_{L^2}\|Au_t \|_{L^2}
		\nonumber\\
		\leq&\
		C\|D\zeta^{\epsilon}\|_{L^{\infty}}\|\nabla u_{t}\|_{L^2}\|\nabla u\|^{\frac{1}{2}}_{L^2}\|Au\|^{\frac{1}{2}}_{L^2}\|Au_{t}\|_{L^2}+\|\zeta^{\epsilon}\|_{L^{\infty}}\|\nabla u_t\|_{L^2}\|Au_{t}\|_{L^2}
		\nonumber\\
		\leq&\ \frac{\nu}{2}\|Au_{t}\|_{L^2}^{2}+C_{\nu}(\|D\zeta^{\epsilon}\|^{2}_{L^{\infty}}\kappa_2 \kappa_4 +\|\zeta^{\epsilon}\|^{2}_{L^{\infty}})\|\nabla u_t \|^{2}_{L^2}.
	\end{align*}
	Then,
	\begin{align}\label{E:3.16}
		\frac{d}{dt}\|\nabla u_t\|^{2}_{L^2}+\nu\|Au_{t}\|_{L^2}^{2} \leq C_{\nu}(\|D\zeta^{\epsilon}\|^{2}_{L^{\infty}}\kappa_2 \kappa_4 +\|\zeta^{\epsilon}\|^{2}_{L^{\infty}})\|\nabla u_t \|^{2}_{L^2}.
	\end{align}
	
	Applying Gronwall inequality to \eqref{E:3.16} over $[s,t+1]$,where $s\in(t,t+1)$, so we can get
	\begin{align*}
		\|\nabla u_t (t+1)\|^2_{L^2}\leq e^{C_{\nu}(\|D\zeta^{\epsilon}\|^{2}_{L^{\infty}}\kappa_2 \kappa_4 +\|\zeta^{\epsilon}\|^{2}_{L^{\infty}})}\|\nabla u_t (s)\|^2_{L^2}.
	\end{align*}
	Furthermore, integrating with respect to $s$ over $[t,t+1]$,
	\begin{align*}
		\|\nabla u_t(t)\|^2_{L^2} \leq e^{C_{\nu}(\|D\zeta^{\epsilon}\|^{2}_{L^{\infty}}\kappa_2 \kappa_4 +\|\zeta^{\epsilon}\|^{2}_{L^{\infty}})}I_4:=\kappa^2_5 ,\; \textup{for} \; \forall t\ge t_3 :=t_2 +1.
	\end{align*}
	
	This completes the proof.
	\end{proof}
	
\begin{lemma}\label{L:3.7}
	For $t\ge 0$, $S(t):V\mapsto V$ is a Lipschitz continuous map on $V$.
\end{lemma}

\begin{proof}
	Let $u_1$ and $u_2$ be two solutions in $V$ with initial values $u_{1,0}$ and $u_{2,0}$, respectively. Set $u=u_1-u_2$, $u_0=u_{1,0}-u_{2,0}$, so
	\begin{align}\label{E:3.17}
		\frac{du}{dt}+\nu Au=B(\zeta^{\epsilon}(u_2)-\zeta^{\epsilon}(u_1),u_2)-B(\zeta^{\epsilon}(u_1),u).
	\end{align}
	Then, we take the inner product of \eqref{E:3.17} with $Au$,
	\begin{align*}
		\frac{1}{2}\frac{d}{dt}\|\nabla u\|^{2}_{L^2}+\nu\|Au\|_{L^2}^{2}=&b(\zeta^{\epsilon}(u_2)-\zeta^{\epsilon}(u_1),u_2,Au)-b(\zeta^{\epsilon}(u_1),u,Au)
		\nonumber\\
		\leq&\
		\|u\|_{L^6}\|\nabla u_2\|_{L^3}\|Au \|_{L^2}+\|\zeta^{\epsilon}(u_1)\|_{L^{\infty}}\|\nabla u \|_{L^2}\|Au \|_{L^2}
		\nonumber\\
		\leq&\
		C\|\nabla u\|_{L^2}\|\nabla u_2\|^{\frac{1}{2}}_{L^2}\|Au_2\|^{\frac{1}{2}}_{L^2}\|Au\|_{L^2}+\|\zeta^{\epsilon}\|_{L^{\infty}}\|\nabla u\|_{L^2}\|Au\|_{L^2}
		\nonumber\\
		\leq&\ \frac{\nu}{2}\|Au\|_{L^2}^{2}+C_{\nu}(\|\nabla u_2\|_{L^{2}}\|Au_2\|_{L^{2}} +\|\zeta^{\epsilon}\|^{2}_{L^{\infty}})\|\nabla u \|^{2}_{L^2}.
	\end{align*}
	So,
	\begin{align}\label{E:3.18}
		\frac{d}{dt}\|\nabla u\|^{2}_{L^2}+\nu\|Au\|_{L^2}^{2} \leq C_{\nu}(\|\nabla u_2\|_{L^{2}}\|Au_2\|_{L^{2}} +\|\zeta^{\epsilon}\|^{2}_{L^{\infty}})\|\nabla u \|^{2}_{L^2},
	\end{align}
	where each term in the parentheses on the right-hand side is integrable on $[0,T]$.
	
	Using Gronwall inequality on \eqref{E:3.18},
	\begin{align}\label{E:3.19}
		\|\nabla u(t)\|^2_{L^2} \leq e^{C_{\nu}\int^{t}_{0}(\|\nabla u_2(s)\|_{L^2}\|Au_2(s)\|_{L^2} +\|\zeta^{\epsilon}\|^{2}_{L^{\infty}})ds}\|\nabla u_0\|^2_{L^2}.
	\end{align}
	
	This completes the proof.
\end{proof}

\section{$(V,H^2)$-Global attractors}

In this section, we recall the definition of global attractors. Simultaneously, we prove that $\{S(t)\}_{t\ge0}$ has a $(V,H^2)$-global attractor.
\begin{definition}\label{D:4.1}
	Assume that $\mathcal{A}$ is a subset of $H^2$. So $\mathcal{A}$ is said to be a $(V,H^2)$-global attractor if the following conditions are satisfied:
	\begin{enumerate}[(i)]
		\item $\mathcal{A}$ is compact in $H^2$.
		\item $\mathcal{A}$ is invariant, i.e., $S(t)\mathcal{A}=\mathcal{A}$, $\forall t \ge 0$.
		\item $\mathcal{A}$ attracts every bounded subset of $V$ with respect to the norm of $H^2$, namely, if $B$ is a bounded set in $V$, then
		\begin{align*}
			dist_{H^2}(S(t)B,\mathcal{A})\rightarrow 0 \; \textup{as} \; t\rightarrow \infty.
		\end{align*}
	\end{enumerate}
\end{definition}

Next, we will introduce some relevant results needed for the proof of the existence of a global attractor for \eqref{E:choudeng} in $H^2$.
\begin{theorem}\label{Th:4.2}
	Assume that $\mathcal{A}$ is a $(V,V)$-global attractor for $\{S(t)\}_{t\ge 0}$, then $\mathcal{A}$ is a $(V,H^2)$-global attractor if and only if
	\begin{enumerate}[(i)]
		\item $\{S(t)\}_{t\ge 0}$ exists a bounded $(V,H^2)$-absorbing set.
		\item $\{S(t)\}_{t\ge 0}$ is $(V,H^2)$-asymptotically compact.
	\end{enumerate}
\end{theorem}

First, we prove that $\{S(t)\}_{t\ge0}$ has a $(V,V)$-global attractor. Then, we show that this attractor is also a $(V,H^2)$-global attractor.

We set
\begin{align*}
	B_1 =\{u\in V:\|\nabla u(t)\|_{L^2}\leq \kappa_2\},
\end{align*}
and
\begin{align*}
	B_2 =\{u\in \mathcal{D}(A):\|A u(t)\|_{L^2}\leq \kappa_4\}.
\end{align*}
By Lemma \ref{L:3.2} and Lemma \ref{L:3.5}, $B_1$ and $B_2$ are bounded $(V,V)$-absorbing set and $(V,H^2)$-absorbing set for $\{S(t)\}_{t\ge0}$, respectively. Furthermore, by the compact embedding $\mathcal{D}(A) \hookrightarrow V$ and Lemma \ref{L:3.5}, we can know that $\{S(t)\}_{t\ge0}$ is $(V,V)$-asymptotically compact. Hence, by the standard attractor theory (see \cite{Temam1997}), $\{S(t)\}_{t\ge0}$ exists a $(V,V)$-global attractor. Moreover, to show that $\mathcal{A}$ is a $(V,H^2)$-global attractor, it suffices to prove that $\{S(t)\}_{t\ge0}$ is $(V,H^2)$-asymptotically compact.

\begin{lemma}\label{L:4.3}
	Assume that $u_{0}\in H$, $f\in H$ and $\zeta^{\epsilon}$ satisfies the conditions in Definition \ref{D:1.1}. Thus, the $\{S(t)\}_{t\ge0}$ is $(V,H^2)$-asymptotically compact.
\end{lemma}

\begin{proof}
	Suppose that $\{u_{0n}\}^{\infty}_{n=1}$ is bounded in $V$, and $t_n \rightarrow \infty$ as $n \rightarrow \infty$. Next, we show that $\{S(t_n)u_{0n}\}^{\infty}_{n=1}$ has a convergent subsequence in $H^2$. Set
	\begin{align*}
		u_n (t)=S(t)u_{0n}, \; w_n(t_n)= \left. \frac{du_n}{dt}\right|_{t=t_n}.
	\end{align*}
	So, \eqref{E:choudeng} is equivalent to the following:
	\begin{align*}
		\nu Au_n(t_n)=f-w_n(t_n)-B(\zeta^{\epsilon}(u_n(t_n)),u_n(t_n)).
	\end{align*}
	By Lemmas \ref{L:3.5} and 3.6, there is $T>0$ such that for $\forall t\ge T$ and $n \in N_{+}$,
	\begin{align}\label{E:4.1}
		\|Au_n(t)\|_{L^2}\le \kappa_4, \; \|\nabla \frac{du_n}{dt}(t)\|_{L^2}\le \kappa_5 ,
	\end{align}
	where the definitions of $\kappa_4$ and $\kappa_5$ are in Lemmas \ref{L:3.5} and \ref{L:3.6}, respectively. Since $t_n \rightarrow \infty$, there exists $N>0$ such that $t_n \ge T$ as $n\ge N$. Therefore, according to \eqref{E:4.1}, for $\forall n\ge N$,
	\begin{align}\label{E:4.2}
		\|Au_n(t_n)\|_{L^2}\le \kappa_4, \; \|\nabla w_n(t_n)\|_{L^2}\le \kappa_5 .
	\end{align}
	Since the compact embeddings $V \hookrightarrow H$, $\mathcal{D}(A) \hookrightarrow V$ and \eqref{E:4.2}, there are $w\in V$ and $v\in \mathcal{D}(A)$ such that, for a subsequence,
	\begin{align}\label{E:4.3}
		w_n(t_n)\rightarrow w \; \textup{strongly in}\; H,
	\end{align}
	\begin{align}\label{E:4.4}
		u_n(t_n)\rightarrow v \; \textup{strongly in}\; V.
	\end{align}
	
	Additionally, we set $\|B(\zeta^{\epsilon}(u_n (t_n)),u_n (t_n))-B(\zeta^{\epsilon}(v),v)\|^2_{L^2}:=R$, so
	\begin{align*}
		R\leq&\
		2\|B(\zeta^{\epsilon}(u_n (t_n)),u_n (t_n)-v)\|^2_{L^2}+2\|B(\zeta^{\epsilon}(u_n (t_n))-\zeta^{\epsilon}(v),v)\|^2_{L^2}
		\nonumber\\
		\leq&\
		2\|\zeta^{\epsilon}(u_n (t_n))\|^2_{L^{\infty}}\|\nabla (u_n (t_n)-v)\|^2_{L^2}+C\|u_n (t_n)-v\|^2_{L^6}\|\nabla v\|_{L^2}\|Av\|_{L^2}
		\nonumber\\
		\leq&\
		C(\|\zeta^{\epsilon}(u_n (t_n))\|^2_{L^{\infty}}+	\|\nabla v\|_{L^2}\|Av\|_{L^2})\|\nabla (u_n (t_n)-v)\|^2_{L^2} \rightarrow 0, \; \textup{as} \; n\rightarrow \infty,
	\end{align*}
	that is,
	\begin{align}\label{E:4.5}
		B(\zeta^{\epsilon}(u_n (t_n)),u_n (t_n))\rightarrow B(\zeta^{\epsilon}(v),v) \; \textup{strongly in} \; H.
	\end{align}
	
	Therefore, by \eqref{E:4.3}, \eqref{E:4.4} and \eqref{E:4.5}, as $n\rightarrow \infty$, we have
	\begin{align*}
		\nu Au_n(t_n)\rightarrow f-w-B(\zeta^{\epsilon}(v),v) \; \textup{strongly in} \; H,
	\end{align*}
	which means that $\{u_n(t_n)\}^{\infty}_{n=1}$ converges to $\frac{1}{\nu}A^{-1}(f-w-B(\zeta^{\epsilon}(v),v))$ in $H^2$.
	
	This completes the proof.
\end{proof}

\noindent\textbf{Proof of Theorem \ref{M:1.3}}
	From what has been discussed above, $\{S(t)\}_{t\ge 0}$ has a $(V,V)$-global attractor $\mathcal{A}$. It follows from Lemma \ref{L:3.5} that the set $B_2 =\{u\in \mathcal{D}(A):\|Au\|_{L^2} \le \kappa_4\}$ is a bounded $(V,H^2)$-absorbing set for $\{S(t)\}_{t\ge0}$. Moreover, Lemma \ref{L:4.3} implies that $\{S(t)\}_{t\ge0}$ is $(V,H^2)$-asymptotically compact. Hence, by Theorem \ref{Th:4.2}, $\mathcal{A}$ is a $(V,H^2)$-global attractor for $\{S(t)\}_{t\ge0}$.

\section{$(H^2,H^3)$-Global attractors}
	
	Similar to the way of the previous section, the existence of $(H^2,H^3)$-global attractor is established below.
\begin{lemma}\label{L:5.1}
	Assume that $u_{0}\in H$, $f\in V$ and $\zeta^{\epsilon}$ satisfies the conditions in Definition \ref{D:1.1}. Thus, there is a constant $\kappa_6>0$ and a time $t_3=t_2 +1$ such that for all $t\geq t_3$,
	\begin{align}
		\|A^{\frac{3}{2}}u(t)\|_{L^2}\leq \kappa_6.
	\end{align}
\end{lemma}

\begin{proof}
	By \eqref{E:choudeng}, we can get
	\begin{align}\label{E:5.2}
		\nu A^{\frac{3}{2}}u=-A^{\frac{1}{2}}u_t -A^{\frac{1}{2}}B(\zeta^{\epsilon}(u),u)+A^{\frac{1}{2}}f,
	\end{align}
	then applying Minkowski's inequality, Definition \ref{D:1.1}, and Lemmas \ref{L:3.2}, \ref{L:3.5} and \ref{L:3.6} to \eqref{E:5.2}, where $t\ge t_3$,
	\begin{align*}
		\nu\|A^{\frac{3}{2}}u\|_{L^2}\leq&\|\nabla u_t \|_{L^2}+\|D\zeta^{\epsilon}(u)(\nabla u)^2\|_{L^2}+\|\zeta^{\epsilon}(u)\Delta u\|_{L^2}+\|\nabla f\|_{L^2}
		\nonumber\\
		\leq&\
		\|\nabla u_t \|_{L^2}+\|D\zeta^{\epsilon}(u)\|_{L^{\infty}}\|\nabla u\|^2_{L^4}+\|\zeta^{\epsilon}(u)\|_{L^2}\|\Delta u\|_{L^2}+\|\nabla f\|_{L^2}
		\nonumber\\
		\leq&\
		\|\nabla u_t \|_{L^2}+C\|D\zeta^{\epsilon}\|_{L^{\infty}}(\|Au\|^{\frac{3}{4}}_{L^2}\|\nabla u\|^{\frac{1}{4}}_{L^2})^2 +\|\zeta^{\epsilon}\|_{L^{\infty}}\|\Delta u\|_{L^2}+\|\nabla f\|_{L^2}
		\nonumber\\
		\leq&\
	\kappa_5 +C\|D\zeta^{\epsilon}\|_{L^{\infty}}\kappa^{\frac{3}{2}}_4 \kappa^{\frac{1}{2}}_2 +\|\zeta^{\epsilon}\|_{L^{\infty}}\kappa_4 +\|\nabla f\|_{L^2}:=\nu\kappa_6.
	\end{align*}
	
	This completes the proof.
\end{proof}

\begin{lemma}\label{L:5.2}
	Assume that $u_{0}\in H$, $f\in V$ and $\zeta^{\epsilon}$ satisfies the conditions in Definition \ref{D:1.1}. Thus, there are constants $\kappa_7>0$, $I_5>0$, and a time $t_4=t_3 +1$ such that for all $t\geq t_4$,
	\begin{align}
		\|A u_{t}(t)\|_{L^2}\leq \kappa_7,
	\end{align}
	\begin{align}
		\int^{t+1}_{t}\|A u_{t}(s)\|_{L^2}^{2}ds\leq I_5 .
	\end{align}
\end{lemma}

\begin{proof}
	Integrating  \eqref{E:3.16} from $t$ to $t+1$ and by Lemma \ref{L:3.6}, where $t\ge t_3$,
	\begin{align*}
		\nu \int^{t+1}_{t}\|A u_{t}(s)\|_{L^2}^{2}ds\leq&\|\nabla u_{t}(t)\|_{L^2}^{2}+C_{\nu}(\|D\zeta^{\epsilon}\|^{2}_{L^{\infty}}\kappa_2 \kappa_4 +\|\zeta^{\epsilon}\|^{2}_{L^{\infty}})\int^{t+1}_{t}\|\nabla u_{t}(s)\|_{L^2}^{2}ds
		\nonumber\\
		\leq&\
		(C_{\nu}\|D\zeta^{\epsilon}\|^{2}_{L^{\infty}}\kappa_2 \kappa_4 +C_{\nu}\|\zeta^{\epsilon}\|^{2}_{L^{\infty}}+1)\kappa_5^2 :=\nu I_5 .
	\end{align*}
	
	Moreover, differentiating \eqref{E:5.2} with respect to $t$ yields
	\begin{align}\label{E:5.5}
		\frac{d}{dt}A^{\frac{1}{2}}u_t + \nu A^{\frac{3}{2}}u_t = -A^{\frac{1}{2}}B(D\zeta^{\epsilon}(u)u_t ,u)-A^{\frac{1}{2}}B(\zeta^{\epsilon}(u),u_t ),
	\end{align}
	
	next, taking the inner product of \eqref{E:5.5} with $A^{\frac{3}{2}}u_t$, we have
	\begin{align*}
		\frac{1}{2}\frac{d}{dt}\|A u_t\|^{2}_{L^2}+\nu\|A^{\frac{3}{2}}u_{t}\|_{L^2}^{2}=M,
	\end{align*}
	where the nonlinearity $M$ is given as follows,
	\begin{align*}
		M=&-b(D^2\zeta^{\epsilon}(u)\nabla u u_t,u,A^{\frac{3}{2}}u_t)-b(D\zeta^{\epsilon}(u)\nabla u_t,u,A^{\frac{3}{2}}u_t)-b(D\zeta^{\epsilon}(u)u_t,\nabla u,A^{\frac{3}{2}}u_t)
		\\
		&-b(D\zeta^{\epsilon}(u)\nabla u ,u_t ,A^{\frac{3}{2}}u_t)-b(\zeta^{\epsilon}(u),\nabla u_t ,A^{\frac{3}{2}}u_t)
		\\
		:=& M_1 +M_2 +M_3 +M_4 +M_5 .
	\end{align*}
	For $M_1$, using H\"older and Gagliardo-Nirenberg inequalities, and Definition \ref{D:1.1}, where $t\ge t_2$,
	\begin{align*}
		|M_1|\leq& \int_{\mathbb{T}^3}|D^2\zeta^{\epsilon}(u)||\nabla u||u_t||\nabla u||A^{\frac{3}{2}}u_t|dx
		\\
		\leq&\|D^2\zeta^{\epsilon}(u)\|_{L^{\infty}}\|\nabla u\|^2_{L^{4}}\|u_t\|_{L^{\infty}}\|A^{\frac{3}{2}}u_t\|_{L^2}
		\\
		\leq&C\|D^2\zeta^{\epsilon}\|_{L^{\infty}}\|\nabla u\|^2_{L^{4}}\|\Delta u_t\|_{L^{2}}\|A^{\frac{3}{2}}u_t\|_{L^2}
		\\
		\leq&C\|D^2\zeta^{\epsilon}\|_{L^{\infty}}\|Au\|^{\frac{3}{2}}_{L^2}\|\nabla u\|^{\frac{1}{2}}_{L^{2}}\|A u_t\|_{L^{2}}\|A^{\frac{3}{2}}u_t\|_{L^2}
		\\
		\leq&\frac{\nu}{10}\|A^{\frac{3}{2}}u_t\|^2_{L^2}+C_{\nu}\|D^2\zeta^{\epsilon}\|^2_{L^{\infty}}\kappa^3_4 \kappa_2 \|A u_t\|^2_{L^{2}}.
	\end{align*}
	The other terms can be obtained similarly, where $t\ge t_2$,
	\begin{align*}
		|M_2|\leq& \int_{\mathbb{T}^3}|D\zeta^{\epsilon}(u)||\nabla u_t||\nabla u||A^{\frac{3}{2}}u_t|dx
		\\
		\leq&\|D\zeta^{\epsilon}(u)\|_{L^{\infty}}\|\nabla u_t\|_{L^{6}}\|\nabla u\|_{L^{3}}\|A^{\frac{3}{2}}u_t\|_{L^2}
		\\
		\leq&C\|D\zeta^{\epsilon}\|_{L^{\infty}}\|A u_t\|_{L^{2}}\|\nabla u\|_{L^{3}}\|A^{\frac{3}{2}}u_t\|_{L^2}
		\\
		\leq&C\|D\zeta^{\epsilon}\|_{L^{\infty}}\|Au_t\|_{L^2}\|Au\|^{\frac{1}{2}}_{L^2}\|\nabla u\|^{\frac{1}{2}}_{L^{2}}\|A^{\frac{3}{2}}u_t\|_{L^2}
		\\
		\leq&\frac{\nu}{10}\|A^{\frac{3}{2}}u_t\|^2_{L^2}+C_{\nu}\|D\zeta^{\epsilon}\|^2_{L^{\infty}}\kappa_4 \kappa_2 \|A u_t\|^2_{L^{2}}.
		\\
		\\
		|M_3|\leq& \int_{\mathbb{T}^3}|D\zeta^{\epsilon}(u)|| u_t||\Delta u||A^{\frac{3}{2}}u_t|dx
		\\
		\leq&\|D\zeta^{\epsilon}(u)\|_{L^{\infty}}\| u_t\|_{L^{\infty}}\|\Delta u\|_{L^{2}}\|A^{\frac{3}{2}}u_t\|_{L^2}
		\\
		\leq&C\|D\zeta^{\epsilon}\|_{L^{\infty}}\|A u_t\|_{L^{2}}\|A u\|_{L^{2}}\|A^{\frac{3}{2}}u_t\|_{L^2}
		\\
		\leq&\frac{\nu}{10}\|A^{\frac{3}{2}}u_t\|^2_{L^2}+C_{\nu}\|D\zeta^{\epsilon}\|^2_{L^{\infty}}\kappa^2_4 \|A u_t\|^2_{L^{2}}.
		\\
		\\
		|M_4|\leq& \int_{\mathbb{T}^3}|D\zeta^{\epsilon}(u)||\nabla u||\nabla u_t||A^{\frac{3}{2}}u_t|dx
		\\
		\leq&\|D\zeta^{\epsilon}(u)\|_{L^{\infty}}\|\nabla u\|_{L^{3}}\|\nabla u_t\|_{L^{6}}\|A^{\frac{3}{2}}u_t\|_{L^2}
		\\
		\leq&C\|D\zeta^{\epsilon}\|_{L^{\infty}}\|A u_t\|_{L^{2}}\|\nabla u\|_{L^{3}}\|A^{\frac{3}{2}}u_t\|_{L^2}
		\\
		\leq&C\|D\zeta^{\epsilon}\|_{L^{\infty}}\|Au_t\|_{L^2}\|Au\|^{\frac{1}{2}}_{L^2}\|\nabla u\|^{\frac{1}{2}}_{L^{2}}\|A^{\frac{3}{2}}u_t\|_{L^2}
		\\
		\leq&\frac{\nu}{10}\|A^{\frac{3}{2}}u_t\|^2_{L^2}+C_{\nu}\|D\zeta^{\epsilon}\|^2_{L^{\infty}}\kappa_4 \kappa_2 \|A u_t\|^2_{L^{2}}.
		\\
		\\
		|M_5|\leq& \int_{\mathbb{T}^3}|\zeta^{\epsilon}(u)||\Delta 	u_t||A^{\frac{3}{2}}u_t|dx
		\\
		\leq&\|\zeta^{\epsilon}\|_{L^{\infty}}\|A u_t\|_{L^{2}}\|A^{\frac{3}{2}}u_t\|_{L^2}
		\\
		\leq&\frac{\nu}{10}\|A^{\frac{3}{2}}u_t\|^2_{L^2}+C_{\nu}\|\zeta^{\epsilon}\|^2_{L^{\infty}} \|A u_t\|^2_{L^{2}}.
    \end{align*}
    After rearranging all terms, which yields
	\begin{align}\label{E:5.6}
		\frac{d}{dt}\|Au_{t}\|_{L^2}^{2} +\nu\|A^{\frac{3}{2}}u_t\|^2_{L^2}\leq K_1 \|Au_t\|^2_{L^2},
	\end{align}
	where
	\begin{align*}
		K_1 :=C_{\nu}(\|D^2\zeta^{\epsilon}\|^2_{L^{\infty}}\kappa^3_4 \kappa_2 +\|D\zeta^{\epsilon}\|^2_{L^{\infty}}\kappa_4 \kappa_2+\|D\zeta^{\epsilon}\|^2_{L^{\infty}}\kappa^2_4 +\|\zeta^{\epsilon}\|^2_{L^{\infty}}),
	\end{align*}
	and it is integrable on $[0,T]$.
	
	Moreover, applying Gronwall inequality to \eqref{E:5.6} over $[s,t+1]$, where $s\in (t,t+1)$,
	\begin{align*}
		\|Au_t (t+1)\|^2_{L^2}\leq e^{K_1}\|Au_t (s)\|^2_{L^2}.
	\end{align*}
	Then, integrating once more with respect to $s$ over $[t,t+1]$, and for any $t\ge t_3 +1 :=t_4$, we have
	\begin{align*}
		\|A u_t(t)\|^2_{L^2} \leq e^{C}I_5:=\kappa^2_7 .
	\end{align*}
	
	This completes the proof.
\end{proof}

\begin{lemma}\label{L:5.3}
	For $t\ge 0$, $S(t):H^2 \mapsto H^2$ is a Lipschitz continuous map on $H^2$.
\end{lemma}

\begin{proof}
	We still let $u_1$ and $u_2$ be two solutions in $H^2$ with initial values $u_{1,0}$ and $u_{2,0}$ separately, and let $u_1 > u_2$. Set $u=u_1-u_2$, $u_0=u_{1,0}-u_{2,0}$, so
	\begin{align}\label{E:5.7}
		\frac{du}{dt}+\nu Au=B(\zeta^{\epsilon}(u_2)-\zeta^{\epsilon}(u_1),u_2)-B(\zeta^{\epsilon}(u_1),u).
	\end{align}
	Then, we take the inner product of \eqref{E:5.7} with $A^2 u$, and take $v\in (u_2 ,u_1)$,
	\begin{align*}
		\frac{1}{2}\frac{d}{dt}\|A u\|^{2}_{L^2}+\nu\|A^{\frac{3}{2}}u\|_{L^2}^{2}=N,
	\end{align*}
	where the nonlinearity is defined as
	\begin{align*}
		N=&-b(D^2 \zeta^{\epsilon}(v)\nabla vu,u_2 ,A^{\frac{3}{2}}u)-b(D \zeta^{\epsilon}(v)\nabla u,u_2 ,A^{\frac{3}{2}}u)+b( \zeta^{\epsilon}(u_2)-\zeta^{\epsilon}(u_1),\nabla u_2 ,A^{\frac{3}{2}}u)
		\\
		&+b(D\zeta^{\epsilon}(u_1)\nabla u_1 ,u ,A^{\frac{3}{2}}u)+b(\zeta^{\epsilon}(u_1),\nabla u,A^{\frac{3}{2}}u)
		\\
		:=&N_1 +N_2 +N_3 +N_4 +N_5 .
	\end{align*}
	For $N_i$, $i=1,2,\cdots,5$, we use the similar methods to estimate them separately,
	\begin{align*}
		|N_1|\leq& \int_{\mathbb{T}^3}|D^2\zeta^{\epsilon}(v)||\nabla v||u||\nabla u_2||A^{\frac{3}{2}}u|dx
		\\
		\leq&\|D^2\zeta^{\epsilon}(v)\|_{L^{\infty}}\|\nabla v\|_{L^{\infty}}\|u\|_{L^{\infty}}\|\nabla u_2\|_{L^2}\|A^{\frac{3}{2}}u\|_{L^2}
		\\
		\leq&C\|D^2\zeta^{\epsilon}\|_{L^{\infty}}\|\nabla v\|_{L^{\infty}}\|Au\|_{L^2}\|\nabla u_2\|_{L^2}\|A^{\frac{3}{2}}u\|_{L^2}
		\\
		\leq&\frac{\nu}{10}\|A^{\frac{3}{2}}u\|^2_{L^2}+C_{\nu}\|D^2\zeta^{\epsilon}\|^2_{L^{\infty}}\|\nabla v\|^2_{L^{\infty}} \|\nabla u_2\|^2_{L^2}\|Au\|^2_{L^2}
		\\
		\leq&\frac{\nu}{10}\|A^{\frac{3}{2}}u\|^2_{L^2}+C_{\nu}\|D^2\zeta^{\epsilon}\|^2_{L^{\infty}}(\|A^{\frac{3}{2}}u_1\|^{\frac{1}{2}}_{L^2}\|Au_1\|^{\frac{1}{2}}_{L^2}+\|A^{\frac{3}{2}}u_2\|^{\frac{1}{2}}_{L^2}\|Au_2\|^{\frac{1}{2}}_{L^2})^2 \|\nabla u_2\|^2_{L^2}\|Au\|^2_{L^2}.
		\end{align*}
		\begin{align*}
		|N_2|\leq& \int_{\mathbb{T}^3}|D\zeta^{\epsilon}(v)||\nabla u||\nabla u_2||A^{\frac{3}{2}}u|dx
		\\
		\leq&\|D\zeta^{\epsilon}(v)\|_{L^{\infty}}\|\nabla u\|_{L^6}\|\nabla u_2\|_{L^{3}}\|A^{\frac{3}{2}}u\|_{L^2}
		\\
		\leq&C\|D\zeta^{\epsilon}\|_{L^{\infty}}\|A u\|_{L^2}\|\nabla u_2 \|_{L^3}\|A^{\frac{3}{2}}u\|_{L^2}
		\\
		\leq&C\|D\zeta^{\epsilon}\|_{L^{\infty}}\|Au_2 \|^{\frac{1}{2}}_{L^2}\|\nabla u_2 \|^{\frac{1}{2}}_{L^2}\|A u\|_{L^2}\|A^{\frac{3}{2}}u\|_{L^2}
		\\
		\leq&\frac{\nu}{10}\|A^{\frac{3}{2}}u\|^2_{L^2}+C_{\nu}\|D\zeta^{\epsilon}\|^2_{L^{\infty}}\|Au_2\|_{L^2}\|\nabla u_2 \|_{L^2}\|Au\|^2_{L^2}.
		\\
		\\
		|N_3|\leq& \int_{\mathbb{T}^3}|\zeta^{\epsilon}(u_2)-\zeta^{\epsilon}(u_1)||\Delta u_2||A^{\frac{3}{2}}u|dx
		\\
		\leq&\|u\|_{L^{\infty}}\|\Delta u_2\|_{L^2}\|A^{\frac{3}{2}}u\|_{L^2}
		\\
		\leq&C\|A u\|_{L^2}\|Au_2 \|_{L^2}\|A^{\frac{3}{2}}u\|_{L^2}
		\\
		\leq&\frac{\nu}{10}\|A^{\frac{3}{2}}u\|^2_{L^2}+C_{\nu}\|Au_2\|^2_{L^2}\|Au\|^2_{L^2}.
		\\
		\\
		|N_4|\leq& \int_{\mathbb{T}^3}|D\zeta^{\epsilon}(u_1)||\nabla u_1||\nabla u||A^{\frac{3}{2}}u|dx
		\\
		\leq&\|D\zeta^{\epsilon}(u_1)\|_{L^{\infty}}\|\nabla u_1\|_{L^3}\|\nabla u\|_{L^6}\|A^{\frac{3}{2}}u\|_{L^2}
		\\
		\leq&C\|D\zeta^{\epsilon}\|_{L^{\infty}}\|\nabla u_1\|_{L^3}\|A u \|_{L^2}\|A^{\frac{3}{2}}u\|_{L^2}
		\\
		\leq&C\|D\zeta^{\epsilon}\|_{L^{\infty}}\|Au_1 \|^{\frac{1}{2}}_{L^2}\|\nabla u_1\|^{\frac{1}{2}}_{L^2}\|A u\|_{L^2}\|A^{\frac{3}{2}}u\|_{L^2}
		\\
		\leq&\frac{\nu}{10}\|A^{\frac{3}{2}}u\|^2_{L^2}+C_{\nu}\|D\zeta^{\epsilon}\|^2_{L^{\infty}}\|Au_1\|_{L^2}\|\nabla u_1\|_{L^2}\|Au\|^2_{L^2}.
		\\
		\\
		|N_5|\leq& \int_{\mathbb{T}^3}|\zeta^{\epsilon}(u_1)||\Delta u||A^{\frac{3}{2}}u|dx
		\\
		\leq&\|\zeta^{\epsilon}(u_1)\|_{L^{\infty}}\|\Delta u\|_{L^2}\|A^{\frac{3}{2}}u\|_{L^2}
		\\
		\leq&C\|\zeta^{\epsilon}\|_{L^{\infty}}\|A u\|_{L^2}\|A^{\frac{3}{2}}u\|_{L^2}
		\\
		\leq&\frac{\nu}{10}\|A^{\frac{3}{2}}u\|^2_{L^2}+C_{\nu}\|\zeta^{\epsilon}\|^2_{L^{\infty}}\|Au\|^2_{L^2}.
	\end{align*}
	So,
	\begin{align}\label{E:5.8}
		\frac{d}{dt}\|Au\|^2_{L^2} +\nu\|A^{\frac{3}{2}}u\|^2_{L^2}\leq K_2(t)\|Au\|^2_{L^2},
	\end{align}
	where
	\begin{align*}
		K_2(t):=&C_{\nu}\big[\|D^2\zeta^{\epsilon}\|^2_{L^{\infty}}(\|A^{\frac{3}{2}}u_1\|^{\frac{1}{2}}_{L^2}\|Au_1\|^{\frac{1}{2}}_{L^2}+\|A^{\frac{3}{2}}u_2\|^{\frac{1}{2}}_{L^2}\|Au_2\|^{\frac{1}{2}}_{L^2})^2 \|\nabla u_2\|^2_{L^2}
		\\
		+&\|D\zeta^{\epsilon}\|^2_{L^{\infty}}\|Au_2\|_{L^2}\|\nabla u_2 \|_{L^2}+\|Au_2\|^2_{L^2}+\|D\zeta^{\epsilon}\|^2_{L^{\infty}}\|Au_1\|_{L^2}\|\nabla u_1\|_{L^2}+\|\zeta^{\epsilon}\|^2_{L^{\infty}}\big],
	\end{align*}
	and $K_2 (t)$ is also integrable on $[0,T]$.
	
	Applying Gronwall inequality to \eqref{E:5.8}, so we have
	\begin{align*}
		\|A u(t)\|^2_{L^2} \leq e^{\int^{t}_{0}K_2(s)ds}\|A u_0\|^2_{L^2}.
	\end{align*}
	
	This completes the proof.
\end{proof}

\begin{definition}\label{D:5.4}
	Assume that subset $\mathcal{A} \subset H^3$, so we say that $\mathcal{A}$ is an $(H^2,H^3)$-global attractor if the following conditions are satisfied:
	\begin{enumerate}[(i)]
		\item $\mathcal{A}$ is compact in $H^3$.
		\item $\mathcal{A}$ is invariant, i.e., $\forall t \ge 0$, $S(t)\mathcal{A}=\mathcal{A}$.
		\item $\mathcal{A}$ attracts all bounded subsets of $H^2$ with respect to the norm of $H^3$, that is, if bounded set $B\subset H^2$, then
		\begin{align*}
			dist_{H^3}(S(t)B,\mathcal{A})\rightarrow 0 \; \textup{as} \; t\rightarrow \infty.
		\end{align*}
	\end{enumerate}
\end{definition}

\begin{theorem}\label{Th:5.5}
	Assume that $\mathcal{A}$ is an $(H^2,H^2)$-global attractor for $\{S(t)\}_{t\ge 0}$, then $\mathcal{A}$ is also an $(H^2,H^3)$-global attractor if and only if
	\begin{enumerate}[(i)]
		\item $\{S(t)\}_{t\ge 0}$ has a bounded $(H^2,H^3)$-absorbing set.
		\item $\{S(t)\}_{t\ge 0}$ is $(H^2,H^3)$-asymptotically compact.
	\end{enumerate}
\end{theorem}

Similar to the proof of the existence of $(V,H^2)$-global attractor, we still show that $\{S(t)\}_{t\ge0}$ exists an $(H^2,H^2)$-global attractor firstly, then we will prove that it also an $(H^2,H^3)$-global attractor by Theorem \ref{Th:5.5}.

$\{S(t)\}_{t\ge0}$ has a $(V,H^2)$-global attractor by Theorem \ref{M:1.3} and $H^2 \hookrightarrow V$, so we can deduce that $\mathcal{A}$ is also an $(H^2,H^2)$-global attractor. Then we set
\begin{align*}
	B_3 =\{u\in \mathcal{D}(A^{\frac{3}{2}}):\|A^{\frac{3}{2}} u(t)\|_{L^2}\leq \kappa_6\}.
\end{align*}
By Lemma \ref{L:5.1}, $B_3$ is a bounded $(H^2,H^3)$-absorbing set for $\{S(t)\}_{t\ge0}$. Since $\mathcal{A}$ is an $(H^2,H^2)$-global attractor, we can know that $\{S(t)\}_{t\ge0}$ is $(H^2,H^2)$-asymptotically compact. Therefore, we only need to show that $\{S(t)\}_{t\ge0}$ is $(H^2,H^3)$-asymptotically compact.

\begin{lemma}\label{L:5.6}
	Assume that $u_{0}\in H$, $f\in V$ and $\zeta^{\epsilon}$ satisfies the conditions in Definition \ref{D:1.1}. Thus, the $\{S(t)\}_{t\ge0}$ is $(H^2,H^3)$-asymptotically compact.
\end{lemma}

\begin{proof}
	We still set that $\{u_{0n}\}^{\infty}_{n=1}$ is bounded in $H^2$, and $t_n \rightarrow \infty$ as $n \rightarrow \infty$. Now, we prove that $\{S(t_n)u_{0n}\}^{\infty}_{n=1}$ exists a convergent subsequence in $H^3$. We set
	\begin{align*}
		u_n (t)=S(t)u_{0n}, \; w_n(t_n)= \left. \frac{du_n}{dt}\right|_{t=t_n}.
	\end{align*}
	So, by \eqref{E:5.2},
	\begin{align*}
		\nu A^{\frac{3}{2}}u_n(t_n)=A^{\frac{1}{2}}f-A^{\frac{1}{2}}w_n(t_n)-A^{\frac{1}{2}}B(\zeta^{\epsilon}(u_n(t_n)),u_n(t_n)).
	\end{align*}
	By Lemmas \ref{L:5.1} and \ref{L:5.2}, there exists $T>0$ such that for all $ t\ge T$ and $n \in N_{+}$,
	\begin{align}\label{E:5.9}
		\|A^{\frac{3}{2}}u_n(t)\|_{L^2}\le \kappa_6, \; \|A\frac{du_n}{dt}(t)\|_{L^2}\le \kappa_7 ,
	\end{align}
	here, the details of $\kappa_6$ and $\kappa_7$ are in Lemmas \ref{L:5.1} and \ref{L:5.2}, respectively. Since $t_n \rightarrow \infty$, there is a $N>0$ such that $t_n \ge T$ as $n\ge N$. So, by \eqref{E:5.9}, for $\forall n\ge N$,
	\begin{align}\label{E:5.10}
		\|A^{\frac{3}{2}}u_n(t_n)\|_{L^2}\le \kappa_6, \; \|Aw_n(t_n)\|_{L^2}\le \kappa_7 .
	\end{align}
	Since \eqref{E:5.10} and the embeddings $\mathcal{D}(A^{\frac{3}{2}}) \hookrightarrow \mathcal{D}(A) \hookrightarrow V$ are compact, there are $v\in \mathcal{D}(A^{\frac{3}{2}})$ and $w\in \mathcal{D}(A)$ such that passing to a subsequence,
	\begin{align}\label{E:5.11}
		u_n(t_n)\rightarrow v \; \textup{strongly in}\; H^2,
	\end{align}
	\begin{align}\label{E:5.12}
		w_n(t_n)\rightarrow w \; \textup{strongly in}\; V.
	\end{align}
	
	Additionally, we set
	\begin{align*}
	E:=&A^{\frac{1}{2}}B(\zeta^{\epsilon}(u_n (t_n)),u_n (t_n))-A^{\frac{1}{2}}B(\zeta^{\epsilon}(v),v)
	\\
	=&B(D\zeta^{\epsilon}(u_n (t_n))\nabla (u_n (t_n)),u_n (t_n)-v)+B(D\zeta^{\epsilon}(u_n (t_n))\nabla (u_n (t_n)-v),v)
	\\
	+&B((D\zeta^{\epsilon}(u_n (t_n))-D\zeta^{\epsilon}(v))\nabla v,v)+B(\zeta^{\epsilon}(u_n (t_n)),\nabla (u_n (t_n)-v))+B(\zeta^{\epsilon}(u_n (t_n))-\zeta^{\epsilon}(v),\nabla v)
	\\
	:=&E_1 +E_2 +E_3 +E_4 +E_5.
	\end{align*}
	Next, taking the $H$-norm of $E_i$($i=1,2,\cdots 5$), and estimating them one by one with the similar methods.
	\begin{align*}
		\|E_1\|_{L^2} \leq& \|D\zeta^{\epsilon}(u_n (t_n))\|_{L^{\infty}}\|\nabla u_n (t_n)\|_{L^3}\|\nabla (u_n (t_n)-v)\|_{L^6}
		\\
		\leq&C\|D\zeta^{\epsilon}\|_{L^{\infty}}\|\nabla u_n (t_n)\|_{L^3}\|A(u_n (t_n)-v)\|_{L^2}
		\\
		\leq&C\|D\zeta^{\epsilon}\|_{L^{\infty}}\|Au_n (t_n)\|^{\frac{1}{2}}_{L^2}\|\nabla u_n (t_n)\|^{\frac{1}{2}}_{L^2}\|A(u_n (t_n)-v)\|_{L^2}.
	\end{align*}
	\begin{align*}
		\|E_2\|_{L^2} \leq& \|D\zeta^{\epsilon}(u_n (t_n))\|_{L^{\infty}}\|\nabla (u_n (t_n)-v)\|_{L^6}\|\nabla v\|_{L^3}
		\\
		\leq&C\|D\zeta^{\epsilon}\|_{L^{\infty}}\|A(u_n (t_n)-v)\|_{L^2}\|\nabla v\|_{L^3}
		\\
		\leq&C\|D\zeta^{\epsilon}\|_{L^{\infty}}\|Av\|^{\frac{1}{2}}_{L^2}\|\nabla v\|^{\frac{1}{2}}_{L^2}\|A(u_n (t_n)-v)\|_{L^2}.
	\end{align*}
	\begin{align*}
		\|E_3\|_{L^2} \leq& \|D^2\zeta^{\epsilon}(m)\|_{L^{\infty}}\|u_n (t_n)-v\|_{L^{\infty}}\|\nabla v\|^{2}_{L^4}
		\\
		\leq&C\|D^2\zeta^{\epsilon}\|_{L^{\infty}}\|A(u_n (t_n)-v)\|_{L^2}\|\nabla v\|^{2}_{L^4}
		\\
		\leq&C\|D^2\zeta^{\epsilon}\|_{L^{\infty}}\|Av\|^{\frac{3}{2}}_{L^2}\|\nabla v\|^{\frac{1}{2}}_{L^2}\|A(u_n (t_n)-v)\|_{L^2},
	\end{align*}
	where $m$ is between $u_n (t_n)$ and $v$.
	\begin{align*}
		\|E_4\|_{L^2} \leq&\|\zeta^{\epsilon}(u_n (t_n))\|_{L^{\infty}}\|A(u_n (t_n)-v)\|_{L^2}
		\\
		\leq&\|\zeta^{\epsilon}\|_{L^{\infty}}\|A(u_n (t_n)-v)\|_{L^2}.
	\end{align*}
	\begin{align*}
		\|E_5\|_{L^2} \leq&\|\zeta^{\epsilon}(u_n (t_n))-\zeta^{\epsilon}(v)\|_{L^{\infty}}\|Av\|_{L^2}
		\\
		\leq&\|u_n (t_n)-v\|_{L^{\infty}}\|Av\|_{L^2}
		\\
		\leq&C\|Av\|_{L^2}\|A(u_n (t_n)-v)\|_{L^2}.
	\end{align*}
	So,
	\begin{align*}
		\|E\|_{L^2}\leq&\|E_1\|_{L^2}+\|E_2\|_{L^2}+\|E_3\|_{L^2}+\|E_4\|_{L^2}+\|E_5\|_{L^2}
		\\
		\leq&K_3\|A(u_n (t_n)-v)\|_{L^2}\rightarrow 0, \; \textup{as} \; n\rightarrow \infty,
	\end{align*}
	where
	\begin{align*}
		K_3:=&C(\|D\zeta^{\epsilon}\|_{L^{\infty}}\|Au_n (t_n)\|^{\frac{1}{2}}_{L^2}\|\nabla u_n (t_n)\|^{\frac{1}{2}}_{L^2}+\|D\zeta^{\epsilon}\|_{L^{\infty}}\|Av\|^{\frac{1}{2}}_{L^2}\|\nabla v\|^{\frac{1}{2}}_{L^2}
		\\
		+&\|D^2\zeta^{\epsilon}\|_{L^{\infty}}\|Av\|^{\frac{3}{2}}_{L^2}\|\nabla v\|^{\frac{1}{2}}_{L^2}+\|\zeta^{\epsilon}\|_{L^{\infty}}+\|Av\|_{L^2}).
	\end{align*}
	Namely,
	\begin{align}\label{E:5.13}
		A^{\frac{1}{2}}B(\zeta^{\epsilon}(u_n (t_n)),u_n (t_n))\rightarrow A^{\frac{1}{2}}B(\zeta^{\epsilon}(v),v) \; \textup{strongly in} \; H.
	\end{align}
	
	Hence, according to \eqref{E:5.11}, \eqref{E:5.12} and \eqref{E:5.13}, as $n\rightarrow \infty$, we can get
	\begin{align*}
		\nu A^{\frac{3}{2}}u_n(t_n)\rightarrow A^{\frac{1}{2}}f-A^{\frac{1}{2}}w-A^{\frac{1}{2}}B(\zeta^{\epsilon}(v),v) \; \textup{strongly in} \; H,
	\end{align*}
	in other words, $\{u_n(t_n)\}^{\infty}_{n=1}$ converges to $\frac{1}{\nu}A^{-\frac{3}{2}}(A^{\frac{1}{2}}f-A^{\frac{1}{2}}w-A^{\frac{1}{2}}B(\zeta^{\epsilon}(v),v))$ in $H^3$.
	
	This completes the proof.
\end{proof}

\noindent\textbf{Proof of Theorem \ref{M:1.4}}
From what has been discussed above, $\{S(t)\}_{t\ge 0}$ has an $(H^2,H^2)$-global attractor $\mathcal{A}$. Since Lemma \ref{L:5.1}, the bounded set $B_3 =\{u\in \mathcal{D}(A^{\frac{3}{2}}):\|A^{\frac{3}{2}}u\|_{L^2} \le \kappa_6\}$ is an $(H^2,H^3)$-absorbing set for $\{S(t)\}_{t\ge0}$. Besides, we deduces that $\{S(t)\}_{t\ge0}$ is $(H^2,H^3)$-asymptotically compact from Lemma \ref{L:5.6}. In summary, because of Theorem \ref{Th:5.5}, $\mathcal{A}$ is also an $(H^2,H^3)$-global attractor for $\{S(t)\}_{t\ge0}$.

\section{Backwards Uniqueness}
In this section, we will deduce that the solutions of equations exist the backwards uniqueness, namely, Theorem \ref{M:1.5} will be confirmed specifically as follows.
\\

\noindent\textbf{Proof of Theorem \ref{M:1.5}}
\begin{proof}
	Let $u_0 ,v_0 \in V$ with $u(t)=S(t)u_0$ and $v(t)=S(t)v_0$, and let $w(t)=u(t)-v(t)$, where $u(t)$ and $v(t)$ are the solutions of \eqref{E:choudeng}. Then, we have $S(T)u_0 =S(T)v_0$ and $w \in L^{\infty}(0,T;V)\cap L^{2}(0,T;\mathcal{D}(A))$. By \eqref{E:choudeng},
	\begin{align*}
		\frac{dw}{dt} + \nu Aw &=-B(\zeta^{\epsilon}(u),w)-B(\zeta^{\epsilon}(u)-\zeta^{\epsilon}(v),v)
		\\
		&:=h(t,w).
	\end{align*}
	And $h(t,w)$ is estimated as follows,
	\begin{align*}
		\|h(t,w)\|_{L^2}&\leq \|B(\zeta^{\epsilon}(u),w)\|_{L^2} +\|B(\zeta^{\epsilon}(u)-\zeta^{\epsilon}(v),v)\|_{L^2}
		\\
		&\leq \|\zeta^{\epsilon}(u)\|_{L^{\infty}}\|\nabla w\|_{L^2} +\|w\|_{L^6}\|\nabla v\|_{L^3}
		\\
		&\leq \|\zeta^{\epsilon}\|_{L^{\infty}}\|\nabla w\|_{L^2} + C\|\nabla w\|_{L^2}\|\nabla v\|^{\frac{1}{2}}_{L^2}\|Av\|^{\frac{1}{2}}_{L^2}
		\\
		&\leq (\|\zeta^{\epsilon}\|_{L^{\infty}}+C\|Av\|_{L^2})\|\nabla w\|_{L^2}
		\\
		&:=k(t)\|\nabla w\|_{L^2},
	\end{align*}
	where $k(t)\in L^{2}(0,T)$.
	
	According to the Appendix A, we notice that $h(t,w)$ satisfies its conditions. Hence, we can get $w(t)=0$, $\forall t\in[0,T]$.
	
	This completes the proof.
\end{proof}

\section{Dimensions of the global attractors}
From the proof of Theorems \ref{M:1.3} and \ref{M:1.4}, we can know that system \eqref{E:choudeng} has a $(V,H^2)$-global attractor $\mathcal{A}$ and an $(H^2,H^3)$-global attractor $\mathcal{A}$. Consequently, $\mathcal{A}$ is also $(V,V)$-global attractor and $(H^2,H^2)$-global attractor, respectively. In this section, we will prove that the global attractor $\mathcal{A}$ satisfies Theorem \ref{M:1.6}.

First, we recall the following theorem (see Theorem 3.2 in \cite{Huang2018} and Lemma 2.6 in \cite{Kou2026}).
\begin{theorem}\label{Th:7.1}
	Assume that $X$ is a Hilbert space, $S(t):X\rightarrow X$ is a mapping, and $\mathcal{A} \subset X$ is an invariant compact set. Assume further that there are constants $l\in [1,+\infty)$ and $\delta \in (0,1)$ such that $\forall a_1, a_2 \in \mathcal{A}$,
	\begin{align*}
		\|S(a_1)-S(a_2)\|_{X}\leq l\|a_1 -a_2\|_{X},
	\\
		\|Q_m [S(a_1)-S(a_2)]\|_{X}\leq \delta \|a_1 -a_2\|_{X},
	\end{align*}
	where $Q_m$ denotes the projection operator from $X$ onto a subspace $(X_m)^{\perp}$ of co-dimension $m\in \mathbb{N}$. Then
	\begin{align*}
		d_{H}(\mathcal{A})\leq d_{F}(\mathcal{A}) \leq m \frac{\ln (\frac{8G^{2}_{a}l^2}{1-\delta^2})}{\ln (\frac{2}{1+\delta^2})},
	\end{align*}
	where $d_{H}(\mathcal{A})$ is the Hausdorff dimension, $d_{F}(\mathcal{A})$ denotes the fractal dimension, and $G_a =\frac{2}{\pi}\int^{1}_{0}\frac{dx}{\sqrt{1-x^4}}=0.8346268\cdots$ is Gauss constant.
\end{theorem}

\noindent\textbf{Proof of Theorem \ref{M:1.6}}
Now, based on the above result, we will use the following theorems to derive Theorem \ref{M:1.6}.
\begin{theorem}\label{Th:7.2}
	On the basis of Theorem \ref{M:1.3}, the global attractor $\mathcal{A}$ exists Hausdorff and fractal dimensions in the $V$ space.
\end{theorem}
\begin{proof}
	Let $u_1(t)=S(t)u^0_1$, $u_2(t)=S(t)u^0_2$ be two solutions of \eqref{E:choudeng}, and $u^0_1 , u^0_2\in \mathcal{A}$. We also set $u(t)=u_1(t)-u_2(t)$, and $u_0 = u^0_1 - u^0_2$, so
	\begin{align}\label{E:7.1}
		\frac{du}{dt} + \nu Au = B(\zeta^{\epsilon}(u_2)-\zeta^{\epsilon}(u_1),u_2)-B(\zeta^{\epsilon}(u_1),u).
	\end{align}
	From \eqref{E:3.19}, we can know
	\begin{align*}
		\|\nabla u(t)\|^2_{L^2} \leq& e^{C_{\nu}\int^{t}_{0}(\|\nabla u_2(s)\|_{L^2}\|Au_2(s)\|_{L^2} +\|\zeta^{\epsilon}\|^{2}_{L^{\infty}})ds}\|\nabla u_0\|^2_{L^2}
		\\
		\leq& e^{C_{\nu}T(\kappa_2 \kappa_4 +\|\zeta^{\epsilon}\|^2_{L^{\infty}})}\|\nabla u_0\|^2_{L^2},
	\end{align*}
	defining $l^2_1:= e^{C_{\nu}T(\kappa_2 \kappa_4 +\|\zeta^{\epsilon}\|^2_{L^{\infty}})}\geq 1$, so we have
	\begin{align}\label{E:7.2}
		\|\nabla u(t)\|^2_{L^2} \leq l^2_1 \|\nabla u_0\|^2_{L^2}.
	\end{align}
	
	Let $P_m$ denotes the orthogonal projection operator from $H$ onto the subspace $span\{\omega_1 , \cdots \omega_m \}$, and $Q_m =I-P_m$, then $u=P_m u+Q_m u$. Applying $Q_m$ to \eqref{E:7.1} and taking the inner product with $AQ_m u$, which yields
	\begin{align*}
		\frac{1}{2}\frac{d}{dt}\|\nabla Q_m u\|^2_{L^2} +\nu \|AQ_m u\|^2_{L^2} =&(B(\zeta^{\epsilon}(u_2)-\zeta^{\epsilon}(u_1),u_2),AQ_m u)-(B(\zeta^{\epsilon}(u_1),u),AQ_m u)
		\\
		\leq& \|u\|_{L^6}\|\nabla u_2\|_{L^3}\|AQ_m u\|_{L^2}+\|\zeta^{\epsilon}(u_1 )\|_{L^{\infty}}\|\nabla u\|_{L^2}\|AQ_m u\|_{L^2}
		\\
		\leq& c\|\nabla u\|_{L^2}\|A u_2 \|_{L^2}\|AQ_m u\|_{L^2}+\|\zeta^{\epsilon}\|_{L^{\infty}}\|\nabla u\|_{L^2}\|AQ_m u\|_{L^2}
		\\
		\leq& \frac{\nu}{2}\|AQ_m u\|^2_{L^2}+C_{\nu}(\kappa_4 +\|\zeta^{\epsilon}\|^2_{L^{\infty}})\|\nabla u\|^2_{L^2},
	\end{align*}
	by the above inequality, and $\lambda_{m+1} \|\nabla Q_m u\|^2_{L^2}\leq \|AQ_m u\|^2_{L^2}$, we can get
	\begin{align}\label{E:7.3}
		\frac{d}{dt}\|\nabla Q_m u\|^2_{L^2} +\nu \lambda_{m+1} \|\nabla Q_m u\|^2_{L^2}\leq C_{\nu}(\kappa_4 +\|\zeta^{\epsilon}\|^2_{L^{\infty}})\|\nabla u\|^2_{L^2}.
	\end{align}
	
	Next, we use the Gronwall inequality to get
	\begin{align}\label{E:7.4}
		\|\nabla Q_m u(T)\|^2_{L^2} \leq \delta^{2}_{m,1} \|\nabla u_0\|^2_{L^2},
	\end{align}
	where $\delta^{2}_{m,1}:= e^{-\nu \lambda_{m+1}T}+C_{\nu}l^2_1 (\kappa_4 +\|\zeta^{\epsilon}\|^2_{L^{\infty}})\frac{1-e^{-\nu \lambda_{m+1}T}}{\nu \lambda_{m+1}}$. And there has a $\delta_1 \in(0,1)$ such that $\delta_{m,1} \leq \delta_1$ when $m$ is sufficiently large.
	
	From \eqref{E:7.2} and \eqref{E:7.4}, it follows that
	\begin{align*}
		\|S(T)u^0_1 -S(T)u^0_2\|_V \leq l_1\|u^0_1 -u^0_2\|_V ,
		\\
		\|Q_m (S(T)u^0_1 -S(T)u^0_2)\|_V \leq \delta_1 \|u^0_1 -u^0_2\|_V .
	\end{align*}
	Hence, from Theorem \ref{Th:7.1}, we can deduce Theorem \ref{Th:7.2} directly.
	
	This completes the proof.
\end{proof}

\begin{theorem}\label{Th:7.3}
	Under the assumptions of Theorem \ref{M:1.4}, the global attractor $\mathcal{A}$ also has Hausdorff and fractal dimensions in the $H^2$ space.
\end{theorem}
\begin{proof}
	We still set that $u_1(t)=S(t)u^0_1$, $u_2(t)=S(t)u^0_2$ be two solutions of \eqref{E:choudeng}, $u^0_1 , u^0_2\in \mathcal{A}$, and $u(t)=u_1(t)-u_2(t)$ with $u_0 = u^0_1 - u^0_2$, then
	\begin{align}\label{E:7.5}
		\frac{du}{dt} + \nu Au = B(\zeta^{\epsilon}(u_2)-\zeta^{\epsilon}(u_1),u_2)-B(\zeta^{\epsilon}(u_1),u),
	\end{align}
	and combined with \eqref{E:5.8}, it can be generated
	\begin{align*}
		\|Au(t)\|^2_{L^2} \leq& e^{\int^t_0 K_2 (s)ds}\|Au_0\|^2_{L^2}
		\\
		\leq& e^{K_{0}T}\|Au_0\|^2_{L^2},
	\end{align*}
	where the constant $K_0 := C_{\nu} (4\kappa^2_2 \kappa_4 \kappa_6 \|D^2\zeta^{\epsilon}\|^2_{L^{\infty}}+2\kappa_2 \kappa_4\|D\zeta^{\epsilon}\|^2_{L^{\infty}}+\kappa^2_4 +\|\zeta^{\epsilon}\|^2_{L^{\infty}})$, and we set $l^2_2 :=e^{K_{0}T} \geq 1$, so
	\begin{align}\label{E:7.6}
		\|Au(t)\|^2_{L^2} \leq l^2_2 \|Au_0\|^2_{L^2}.
	\end{align}
	
	Still taking the $P_m$ and $Q_m =I-P_m$ in Theorem \ref{Th:7.2}, and $u=P_m u+Q_m u$. It also acts $Q_m$ on \eqref{E:7.5}, and then making the inner product with $A^2 Q_m u$, which yields
	\begin{align*}
		\frac{1}{2}\frac{d}{dt}\|AQ_m  u\|^{2}_{L^2}+\nu\|A^{\frac{3}{2}}Q_m u\|_{L^2}^{2}=(B(\zeta^{\epsilon}(u_2)-\zeta^{\epsilon}(u_1),u_2),A^2Q_m u)-(B(\zeta^{\epsilon}(u_1),u),A^2Q_m u).
	\end{align*}
	Similar to the proof process in Lemma \ref{L:5.3}, it can be obtained
	\begin{align*}
		\frac{d}{dt}\|AQ_m  u\|^{2}_{L^2}+\nu\|A^{\frac{3}{2}}Q_m u\|_{L^2}^{2} \leq K_2 (t)\|Au\|^2_{L^2},
	\end{align*}
	where the detail of $K_2 (t)$ is given by Lemma \ref{L:5.3}. And since $\lambda_{m+1} \|A Q_m u\|^2_{L^2}\leq \|A^{\frac{3}{2}}Q_m u\|^2_{L^2}$,
	\begin{align}\label{E:7.7}
		\frac{d}{dt}\|AQ_m  u\|^{2}_{L^2}+\nu \lambda_{m+1}\|AQ_m u\|_{L^2}^{2} \leq K_0 \|Au\|^2_{L^2}.
	\end{align}
	
	Applying Gronwall inequality to \eqref{E:7.7},
	\begin{align*}
		\|AQ_m  u(T)\|^{2}_{L^2}\leq \delta^2_{m,2}\|Au_0\|^{2}_{L^2},
	\end{align*}
	letting also $\delta^2_{m,2}:=e^{-\nu \lambda_{m+1}T}+l^2_2 K_0 \frac{1-e^{-\nu \lambda_{m+1}T}}{\nu \lambda_{m+1}}$. Similarly, when $m$ is large enough, there is a constant $\delta_{2}\in(0,1)$ such that $\delta^2_{m,2}\leq \delta^2_{2}$. Namely,
	\begin{align}\label{E:7.8}
		\|AQ_m  u(T)\|^{2}_{L^2}\leq \delta^2_{2}\|Au_0\|^2_{L^2}.
	\end{align}
	
	Then, by \eqref{E:7.6} and \eqref{E:7.8},
	\begin{align*}
		\|S(T)u^0_1 -S(T)u^0_2\|_{H^2} \leq l_2\|u^0_1 -u^0_2\|_{H^2} ,
		\\
		\|Q_m (S(T)u^0_1 -S(T)u^0_2)\|_{H^2} \leq \delta_2 \|u^0_1 -u^0_2\|_{H^2}.
	\end{align*}
	Therefore, we can also get the conclusion that we need.
	
	This completes the proof.
\end{proof}

\section{Determining nodes and determining modes}
In this section, we will discuss that the conclusion on determining nodes and determining modes, as stated in Theorems \ref{M:1.7} and \ref{M:1.8}. First, we should recall the relevant definitions and preliminary lemmas (see \cite{Robinson2013}), and give the proofs by the previous estimations.
\begin{definition}\label{D:8.1}
	Let $\mathcal{N}=\{x_j\}^N_{j=1}\subset \mathbb{T}^3$ be a set consisting of finitely many nodes. For any two orbits $u_1 (t)$ and $u_2 (t)$ with initial values in the global attractor $\mathcal{A}$, if
	\begin{align*}
		\lim_{t\rightarrow \infty} \sup_{j} |u_1(x_j ,t)-u_2(x_j ,t)|=0
	\end{align*}
	implies
	\begin{align*}
		\lim_{t\rightarrow \infty} \|u_1(x_j ,t)-u_2(x_j ,t)\|_{L^{\infty}}=0,
	\end{align*}
	then $\mathcal{N}$ is called a determining nodes set on $\mathcal{A}$.
\end{definition}
\begin{lemma}\label{L:8.2}(see Lemma 35 in \cite{Robinson2013})
	Setting
	\begin{align*}
		d(\mathcal{N})=\sup_{x\in \mathbb{T}^3}\min_{j}|x-x_j|,
	\end{align*}
	there exists a $x_j$ such that $|x-x_j|\leq d(\mathcal{N})$ for any $x\in \mathbb{T}^3$.
	
	And if $\forall w\in \mathcal{D}(A)$, defining
	\begin{align*}
		\eta(w)=\max_{x_j \mathcal{N}}|w(x_j)|,
	\end{align*}
	then there has a positive constant $C_M$, such that
	\begin{align*}
		\|w\|_{L^{\infty}}\leq \eta(w) +C_M d(\mathcal{N})^{\frac{1}{4}}\|Aw\|_{L^2}.
	\end{align*}
\end{lemma}
\begin{proof}
	This lemma can be derived directly from the Sobolev embedding and the above definitions.
\end{proof}

\noindent\textbf{Proof of Theorem \ref{M:1.7}}

\begin{proof}
	Taking any $u^0_1 , u^0_2 \in \mathcal{A}$, setting $u_1(t)=S(t)u^0_1 $ and $u_2(t)=S(t)u^0_2 $, letting also $u(t)$ be the difference of $u_1(t)$ and $u_2(t)$, so
	\begin{align}\label{E:8.1}
		\frac{du}{dt} + \nu Au = B(\zeta^{\epsilon}(u_2)-\zeta^{\epsilon}(u_1),u_2)-B(\zeta^{\epsilon}(u_1),u).
	\end{align}
	Taking the inner product of the above equation with \eqref{E:8.1},
	\begin{align*}
		\frac{1}{2}\frac{d}{dt}\|\nabla u\|^2_{L^2} +\nu\|Au\|^2_{L^2}=&(B(\zeta^{\epsilon}(u_2)-\zeta^{\epsilon}(u_1),u_2),Au)-(B(\zeta^{\epsilon}(u_1),u),Au)
		\\
		\leq&\|u\|_{L^{\infty}}\|\nabla u_2\|_{L^2}\|Au\|_{L^2}+\frac{1}{\sqrt{\lambda_1}}\|\zeta^{\epsilon}\|_{L^{\infty}}\|Au\|^2_{L^2}
		\\
		\leq&\eta(u)\|\nabla u_2\|_{L^2}\|Au\|_{L^2}+(C_M \kappa_2 \delta^{\frac{1}{4}} +\frac{1}{\sqrt{\lambda_1}}\|\zeta^{\epsilon}\|_{L^{\infty}})\|Au\|^2_{L^2},
	\end{align*}
	so we can obtain
	\begin{align*}
		\frac{1}{2}\frac{d}{dt}\|\nabla u\|^2_{L^2}+\lambda_1 (\nu-C_M \kappa_2 \delta^{\frac{1}{4}}-\frac{1}{\sqrt{\lambda_1}}\|\zeta^{\epsilon}\|_{L^{\infty}})\|\nabla u\|^2_{L^2}\leq \kappa_2 \kappa_4 \eta(u),
	\end{align*}
	now we select $\delta$ such that $\mu:=\lambda_1 (\nu-C_M \kappa_2 \delta^{\frac{1}{4}}-\frac{1}{\sqrt{\lambda_1}}\|\zeta^{\epsilon}\|_{L^{\infty}})>0$, then
	\begin{align*}
		\frac{1}{2}\frac{d}{dt}\|\nabla u\|^2_{L^2}+\mu \|\nabla u\|^2_{L^2}\leq \kappa_2 \kappa_4 \eta(u).
	\end{align*}
	We can know that $\eta(u)\rightarrow 0$ by the assumption, and we have $\|\nabla u\|^2_{L^2}\leq \kappa^2_2$ when $t$ is large enough.
	
	Hence, by Lemma \ref{L:2.1}, we have $\|\nabla u\|^2_{L^2}\rightarrow 0$ when $t\rightarrow \infty$. Since also $\|\nabla u\|_{L^{\infty}}\leq c \|\nabla u\|^{\frac{1}{4}}_{L^2}\|\nabla u\|^{\frac{3}{4}}_{H^2}$ holds and Theorem \ref{M:1.4} gives $\sup\|u(t)\|_{H^3}\leq C$,it follows that $\|\nabla u\|^2_{L^{\infty}}\rightarrow 0$. Moreover, $\mathcal{N}$ is determining.
	
	This completes the proof.
\end{proof}

\noindent\textbf{Proof of Theorem \ref{M:1.8}}

	\begin{proof}
		We set $u=P_n u+Q_n u = p+q$, and form the inner product of \eqref{E:8.1} with $q$,
		\begin{align*}
			\frac{1}{2}\frac{d}{dt}\|q\|^2_{L^2}+\nu \|\nabla q\|^2_{L^2}=&(B(\zeta^{\epsilon}(u_2)-\zeta^{\epsilon}(u_1),u_2),q)-(B(\zeta^{\epsilon}(u_1),p),q)-(B(\zeta^{\epsilon}(u_1),q),q)
			\\
			\leq& \|p\|_{L^3}\|\nabla u_2\|_{L^6}\|q\|_{L^2}+\|q\|^2_{L^4}\|\nabla u_2\|_{L^2}
			\\
			+&\|\zeta^{\epsilon}(u_1)\|_{L^{\infty}}\|\nabla p\|_{L^2}\|q\|_{L^2}+\|\zeta^{\epsilon}(u_1)\|_{L^{\infty}}\|\nabla q\|_{L^2}\|q\|_{L^2}
			\\
			\leq& c\|p\|^{\frac{1}{2}}_{L^2}\|\nabla p\|^{\frac{1}{2}}_{L^2}\|\nabla u_2\|^{\frac{1}{2}}_{L^2}\|A u_2\|^{\frac{1}{2}}_{L^2}\|q\|_{L^2}+c\|q\|^{\frac{1}{2}}_{L^2}\|\nabla q\|^{\frac{3}{2}}_{L^2}\|\nabla u_2\|_{L^2}
			\\
			+&\|\zeta^{\epsilon}\|_{L^{\infty}}\|\nabla p\|_{L^2}\|q\|_{L^2}+\|\zeta^{\epsilon}\|_{L^{\infty}}\|\nabla q\|_{L^2}\|q\|_{L^2}
			\\
			\leq&\|q\|^2_{L^2}+(c\lambda^{-\frac{1}{4}}_{n+1}\kappa_2 +\lambda^{-\frac{1}{2}}_{n+1}\|\zeta^{\epsilon}\|_{L^{\infty}})\|\nabla q\|^2_{L^2}+c\lambda_n (\kappa_2 \kappa_4+\|\zeta^{\epsilon}\|^2_{L^{\infty}})\|p\|^2_{L^2}.
		\end{align*}
		Also, since$\lambda_{n+1}\|q\|^2_{L^2}\leq \|\nabla q\|^2_{L^2}$,
		\begin{align*}
			\frac{1}{2}\frac{d}{dt}\|q\|^2_{L^2}+c_0 \|q\|^2_{L^2}\leq +c\lambda_n (\kappa_2 \kappa_4+\|\zeta^{\epsilon}\|^2_{L^{\infty}})\|p\|^2_{L^2} ,
		\end{align*}
		here, $c_0 := \lambda_{n+1}\nu -c\lambda^{\frac{3}{4}}_{n+1}\kappa_2 -\lambda^{\frac{1}{2}}_{n+1}\|\zeta^{\epsilon}\|_{L^{\infty}}-1$.
		
		If $c_0 >0$ for sufficiently large $n$, then, by assumption, $\|p\|^2_{L^2}\rightarrow 0$ as $t\rightarrow \infty$. And hence we can get $\|q\|_{L^2}\rightarrow 0$ from Lemma \ref{L:2.1}.
		
		This completes the proof.
	\end{proof}

\section{Appendix}
In the following, we shall complete the proof by using the method in Robinson (\cite{Robinson2013}, Chapter 2 Section 5).

Define
\begin{align*}
	\Lambda (w)=\frac{\|\nabla w\|^2_{L^2}}{\|w\|^2_{L^2}},
\end{align*}
then differentiate $\Lambda (w)$ with respect to $t$, and apply the equivalent forms of $(Aw-\Lambda w,w)=0$(see Lemma 6.1 in \cite{Temam1997}), which yields
\begin{align*}
	\frac{d}{dt}\Lambda (w)&=\frac{-2\nu (\|Aw\|^2_{L^2}-\Lambda^2 \|w\|^2_{L^2})\|w\|^2_{L^2}+2(Aw,h)\|w\|^2_{L^2}-2(w,h)\|\nabla w\|^2_{L^2}}{\|w\|^4_{L^2}}
	\\
	&=-\frac{2\nu}{\|w\|^2_{L^2}}\|Aw-\Lambda w\|^2_{L^2}+\frac{2}{\|w\|^2_{L^2}}(Aw-\Lambda w,h)
	\\
	&\leq -\frac{2\nu}{\|w\|^2_{L^2}}\|Aw-\Lambda w\|^2_{L^2}+\frac{2}{\|w\|^2_{L^2}}\|Aw-\Lambda w\|_{L^2}k(t)\|\nabla w\|_{L^2}
	\\
	&\leq -2\nu(\frac{\|Aw-\Lambda w\|_{L^2}}{\|w\|_{L^2}})^2 + 2k(t)\Lambda^{\frac{1}{2}}(\frac{\|Aw-\Lambda w\|_{L^2}}{\|w\|_{L^2}}),
\end{align*}
rearranging the above inequality, we can get
\begin{align}\label{E:A.1}
	\frac{d}{dt}\Lambda\leq \frac{1}{2\nu}k^2(t)\Lambda.
\end{align}
Then applying Gronwall inequality to \eqref{E:A.1},
\begin{align}
	\Lambda(t)\leq e^{\frac{1}{2\nu}\int^t_0 k^2(s)ds}\Lambda (0).
\end{align}

Now, using a method similar to that in Lemma 5 of \cite{Robinson2013}, we give the following lemma.
\begin{lemma}\label{L:A.1}
	If $w \in L^{\infty}(0,T;V)\cap L^{2}(0,T;\mathcal{D}(A))$, and there exists $T>0$ such that $w(T)=0$, then $w(t)=0$, $\forall t\in [0,T]$.
\end{lemma}
\begin{proof}
	Assume that there has $t_0 \in [0,T)$, such that $\|w(t_0)\|_{L^2} \neq 0$. By the continuity of $w(t)$ in $L^2$, there exists $\epsilon >0$ such that  $w(t_0) \neq 0$ on $[t_0 ,t_0 +\epsilon)$. Simultaneously, let $[t_0 ,t_1)$ be the maximal time interval on which $w(t_0) \neq 0$, where $t_1 \leq T$ and $w(t_1)=0$.
	
	Then, we differentiate $\ln \frac{1}{\|w\|_{L^2}}$ with respect to $t$ on $[t_0 ,t_1)$,
	\begin{align*}
		\frac{d}{dt}\ln\frac{1}{\|w\|_{L^2}} &\leq \nu \Lambda(w) +\frac{|(h,w)|}{\|w\|^2_{L^2}}
		\\
		&\leq \nu \Lambda(w) +k(t)\Lambda^{\frac{1}{2}}(w)
		\\
		&\leq 2\nu \Lambda(w) + C_{\nu}k^2 (t),
	\end{align*}
integrating the above inequality over $[t_0 ,t]$, where $t\in [t_0 ,t_1)$, which yields
\begin{align*}
	\ln\frac{1}{\|w\|_{L^2}} \leq \ln\frac{1}{\|w(t_0)\|_{L^2}} +2\nu \Lambda(0)  e^{\frac{1}{2\nu}\int^T_0 k^2(s)ds} + C_{\nu}\int^T_0 k^2 (s) ds.
\end{align*}
This shows that $\ln\frac{1}{\|w\|_{L^2}} $ is still bounded on $t\rightarrow t_1$, contradicting the assumption. Therefore, since $t_0 \in [0,T]$ is arbitrary, we obtain $w(t)\equiv 0$.

This completes the proof.
\end{proof}
\section*{Acknowledgment}
 H. Liu was supported by the National Natural Science Foundation of China (No. 12271293), Natural Science Foundation of Shandong Province (No.  ZR2024MA069) and the project of Youth Innovation Team of Universities of Shandong Province (No. 2023KJ204). C.F. Sun was supported by the the National Natural Science Foundation of China (No. 11701269) and Natural Science Foundation of Jiangsu Province (No.  BK20231301). J. Xin was  supported  by the Natural Science Foundation of Shandong Province (No. ZR2023MA002).

\section*{Data availability}
All data generated or analysed in this study are incorporated in the present article.

\section*{Conflict of interest}
The authors declare that they have no conflict of interest.

\end{document}